\documentclass[a4paper,12pt]{article}
\usepackage[T2A]{fontenc}
\usepackage[cp1251]{inputenc}
\begin{document}

\author{S.V. Ludkovsky.}

\title{Structure of wrap groups of quaternion and octonion fiber bundles.}

\date{25 April 2008}
\maketitle

\begin{abstract}
This article is devoted to the investigation of structure of wrap
groups of connected fiber bundles over the fields of real $\bf R$,
complex $\bf C$ numbers, the quaternion skew field $\bf H$ and the
octonion algebra $\bf O$. Iterated wrap groups are studied as well.
Their smashed products are constructed.
\end{abstract}

\section{Introduction.}
\par Geometric loop groups of circles were first introduced by
Lefshetz in 1930-th and then their construction was reconsidered by
Milnor in 1950-th. Lefshetz has used the $C^0$-uniformity on
families of continuous mappings, which led to the necessity of
combining his construction with the structure of a free group with
the help of words. Later on Milnor has used the Sobolev's
$H^1$-uniformity, that permitted to introduce group structure more
naturally \cite{milmorse}. \par  The construction of Lefshetz is
very restrictive, because it works with the $C^0$ uniformity of
continuous mappings in compact-open topology. Even for spheres $S^n$
of dimension $n>1$ it does not work directly, but uses the iterated
loop group construction of circles. Then their constructions were
generalized for fibers over circles and spheres with parallel
transport structures over $\bf C$. Smooth Deligne cohomologies were
studied on such groups \cite{gaj}.

\par Wrap groups of quaternion and octonion fibers
as well as for wider classes of fibers over $\bf R$ or $\bf C$ were
defined and various examples were given together with basic theorems
in \cite{luwrgfbqo}. In that paper a construction of wrap groups was
done with the help of Sobolev uniformities, that has permitted to
consider wide families of manifolds and fiber bundles. This paper
continues previous works of the author on this theme
\cite{ludan,lugmlg,lujmslg,lufoclg,luwrgfbqo}. Wrap groups are
generalizations of geometric loop groups from spheres onto a wider
class of manifolds and fiber bundles over them.
\par  Geometric loop groups have important applications in modern physical
theories (see \cite{ish,mensk} and references therein). Groups of
loops are also intensively used in gauge theory. Wrap groups can be
used in the membrane theory which is the generalization of the
string (superstring) theory.
\par In paper \cite{luwrgfbqo} wrap groups of fiber bundles
over quaternions and octonions were defined  and investigated and
numerous examples were described. This article is devoted to
investigation of their structure and uses notations and results of
the previous work. Besides quaternion and octonion also real and
complex fiber bundles are considered with wrap groups for them.
Smashed products of wrap groups are constructed. Iterated wrap
groups are studied as well.

\par All main results of this paper are obtained for the first time
and they are given in Theorems 2, 6, 9, 10, 20, 21, Propositions 3,
7, 8, 12, 13, 17 and Corollary 11.

\par Remind the basic definitions and notations.
\par {\bf 1. Note.} Denote by ${\cal A}_r$ the
Cayley-Dickson algebra such that ${\cal A}_0=\bf R$, ${\cal A}_1=\bf
C$, ${\cal A}_2=\bf H$ is the quaternion skew field, ${\cal A}_3=\bf
O$ is the octonion algebra. Henceforth we consider only $0\le r\le
3$.

\par {\bf 2.1. Remark.} If $M$ is a metrizable space and $K=K_M$
is a closed subset in $M$ of codimension $codim_{\bf R} ~N\ge 2$
such that $M\setminus K =M_1$ is a manifold with corners over ${\cal
A}_r$, then we call $M$ a pseudo-manifold over ${\cal A}_r$, where
$K_M$ is a critical subset.
\par Two pseudo-manifolds $B$ and $C$ are called diffeomorphic, if
$B\setminus K_B$ is diffeomorphic with $C\setminus K_C$ as for
manifolds with corners (see also \cite{gaj,michor}).

\par Take on $M$ a Borel $\sigma $-additive measure $\nu $
such that $\nu $ on $M\setminus K$ coincides with the Riemann volume
element and $\nu (K)=0$, since the real shadow of $M_1$ has it.
\par The uniform space $H^t_p(M_1,N)$ of all continuous piecewise
$H^t$ Sobolev mappings from $M_1$ into $N$ is introduced in the
standard way \cite{lujmslg,lufoclg}, which induces $H^t_p(M,N)$ the
uniform space of continuous piecewise $H^t$ Sobolev mappings on $M$,
since $\nu (K)=0$, where ${\bf R}\ni t\ge [m/2]+1$, $m$ denotes the
dimension of $M$ over $\bf R$, $[k]$ denotes the integer part of
$k\in \bf R$, $[k]\le k$. Then put $H^{\infty }_p(M,N)=\bigcap_{t>m}
H^t_p(M,N)$ with the corresponding uniformity.
\par For manifolds over ${\cal A}_r$ with $1\le r\le 3$ take
as $H^t_p(M,N)$ the completion of the family of all continuous
piecewise ${\cal A}_r$-holomorphic mappings from $M$ into $N$
relative to the $H^t_p$ uniformity, where $[m/2]+1\le t\le \infty $.
Henceforth we consider pseudo-manifolds with connecting mappings of
charts continuous in $M$ and $H^{t'}_p$ in $M\setminus K_M$ for
$0\le r\le 3$, where $t'\ge t$.
\par {\bf 2.2. Note.} Since the octonion algebra $\bf O$
is non-associative, we consider a non-associative subgroup $G$ of
the family $Mat_q({\bf O})$ of all square $q\times q$ matrices with
entries in $\bf O$. More generally $G$ is a group which has a
$H^t_p$ manifold  structure over ${\cal A}_r$ and group's operations
are $H^t_p$ mappings. The $G$ may be non-associative for $r=3$, but
$G$ is supposed to be alternative, that is, $(aa)b=a(ab)$ and
$a(a^{-1}b)=b$ for each $a, b\in G$.

\par As a generalization of pseudo-manifolds there is used the
following (over $\bf R$ and $\bf C$ see \cite{gaj,souriau}). Suppose
that $M$ is a Hausdorff topological space of covering dimension $dim
~M=m$ supplied with a family $ \{ h : U\to M \} $ of the so called
plots $h$ which are continuous maps satisfying conditions $(D1-D4)$:
\par $(D1)$ each plot has as a domain a convex subset $U$ in ${\cal
A}_r^n$, $n\in \bf N$;
\par $(D2)$ if $h: U\to M$ is a plot, $V$ is a convex subset in
${\cal A}_r^l$ and $g: V\to U$ is an $H^t_p$ mapping, then $h\circ
g$ is also a plot, where $t\ge [m/2]+1$;
\par $(D3)$ every constant map from a convex set $U$ in ${\cal
A}_r^n$ into $M$ is a plot;
\par $(D4)$ if $U$ is a convex set in ${\cal A}_r^n$ and
$ \{ U_j: j\in J \} $ is a covering of $U$ by convex sets in ${\cal
A}_r^n$, each $U_j$ is open in $U$, $h: U\to M$ is such that each
its restriction $h|_{U_j}$ is a plot, then $h$ is a plot. Then $M$
is called an $H^t_p$-differentiable space.
\par A mapping $f: M\to N$ between two $H^t_p$-differentiable spaces
is called differentiable if it continuous and for each plot $h: U\to
M$ the composition $f\circ h: U\to N$ is a plot of $N$. A
topological group $G$ is called an $H^t_p$-differentiable group if
its group operations are $H^t_p$-differentiable mappings.
\par  Let $E$, $N$, $F$ be $H^{t'}_p$-pseudo-manifolds or
$H^{t'}_p$-differentiable spaces over ${\cal A}_r$, let also $G$ be
an $H^{t'}_p$ group over ${\cal A}_r$, $t\le t'\le \infty $. A fiber
bundle $E(N,F,G,\pi ,\Psi )$ with a fiber space $E$, a base space
$N$, a typical fiber $F$ and a structural group $G$ over ${\cal
A}_r$, a projection $\pi : E\to N$ and an atlas $\Psi $ is defined
in the standard way \cite{gaj,michor,sulwint} with the condition,
that transition functions are of $H^{t'}_p$ class such that for
$r=3$ a structure group may be non-associative, but alternative.

\par Local trivializations $\phi _j\circ \pi \circ \Psi _k^{-1}:
V_k(E)\to V_j(N)$ induce the $H^{t'}_p$-uniformity in the family $W$
of all principal $H^{t'}_p$-fiber bundles $E(N,G,\pi ,\Psi )$, where
$V_k(E) = \Psi _k(U_k(E))\subset X^2(G)$, $V_j(N) = \phi
_j(U_j(N))\subset X(N)$, where $X(G)$ and $X(N)$ are ${\cal
A}_r$-vector spaces on which $G$ and $N$ are modelled, $(U_k(E),\Psi
_k)$ and $(U_j(N),\phi _j)$ are charts of atlases of $E$ and $N$,
$\Psi _k =\Psi _k^{E}$, $\phi _j = \phi _j^{N}$.

\par If $G=F$ and $G$ acts on itself by left shifts, then a fiber bundle is
called the principal fiber bundle and is denoted by $E(N,G,\pi ,\Psi
)$. As a particular case there may be $G={\cal A}_r^*$, where ${\cal
A}_r^*$ denotes the multiplicative group ${\cal A}_r\setminus \{ 0
\} $. If $G=F= \{ e \} $, then $E$ reduces to $N$.

\par {\bf 3. Definitions.} Let $M$ be a connected
$H^t_p$-pseudo-manifold over ${\cal A}_r$, $0\le r\le 3$ satisfying
the following conditions:

\par $(i)$ it is compact;

\par $(ii)$ $M$ is a union of two closed subsets
over ${\cal A}_r$ $A_1$ and $A_2$, which are pseudo-manifolds and
which are canonical closed subsets in $M$ with $A_1\cap A_2=\partial
A_1\cap
\partial A_2=:A_3$ and a codimension over $\bf R$ of $A_3$ in $M$ is
$codim_{\bf R}A_3=1$, also $A_3$ is a pseudo-manifold;

\par $(iii)$ a finite set of marked points $s_{0,1},...,s_{0,k}$
is in $\partial A_1\cap \partial A_2$, moreover, $\partial A_j$ are
arcwise connected $j=1, 2$;

\par $(iv)$ $A_1\setminus \partial A_1$ and $A_2\setminus \partial A_2$
are $H^t_p$-diffeomorphic with $M\setminus [ \{ s_{0,1},...,s_{0,k}
\} \cup (A_3\setminus Int (\partial A_1\cap \partial A_2))]$ by
mappings $F_j(z)$, where $j=1$ or $j=2$, $\infty \ge t\ge [m/2]+1$,
$m=dim_{\bf R}M$ such that $H^t\subset C^0$ due to the Sobolev
embedding theorem \cite{miha}, where the interior $Int (\partial
A_1\cap \partial A_2)$ is taken in $\partial A_1\cup \partial A_2$.

\par Instead of $(iv)$ we consider also the case
\par $(iv')$  $M$, $A_1$ and $A_2$ are such that
$(A_j\setminus \partial A_j)\cup \{ s_{0,1},...,s_{0,k} \} $ are \\
$C^0([0,1],H^t_p(A_j,A_j))$-retractable on $X_{0,q}\cap A_j$, where
$X_{0,q}$ is a closed arcwise connected subset in $M$, $j=1$ or
$j=2$, $s_{0,q}\in X_{0,q}$, $X_{0,q}\subset K_M$, $q=1,...,k$,
$codim_{\bf R}~ K_M\ge 2$.
\par Let $\hat M$ be a compact connected $H^t_p$-pseudo-manifold
which is a canonical closed subset in ${\cal A}_r^l$ with a boundary
$\partial {\hat M}$ and marked points $\{ {\hat s}_{0,q}\in \partial
{\hat M}: q=1,...,2k \} $ and an $H^t_p$-mapping $\Xi : {\hat M}\to
M$ such that \par $(v)$ $\Xi $ is surjective and bijective from
${\hat M}\setminus \partial {\hat M}$ onto $M\setminus \Xi (\partial
{\hat M})$ open in $M$, $\Xi ({\hat s}_{0,q})=\Xi ({\hat s}_{0,k+q})
= s_{0,q}$ for each $q=1,...,k$, also $\partial M\subset \Xi
(\partial {\hat M})$.

\par A parallel transport structure on a $H^{t'}_p$-differentiable
principal $G$-bundle $E(N,G,\pi ,\Psi )$ with arcwise connected $E$
and $G$ for $H^t_p$-pseudo-manifolds $M$ and $\hat M$ as above over
the same ${\cal A}_r$ with $t'\ge t+1$ assigns to each $H^t_p$
mapping $\gamma $ from $M$ into $N$ and points $u_1,...,u_k\in
E_{y_0}$, where $y_0$ is a marked point in $N$, $y_0=\gamma
(s_{0,q})$, $q=1,...,k$, a unique $H^t_p$ mapping ${\bf P}_{{\hat
\gamma },u}: {\hat M}\to E$ satisfying conditions $(P1-P5)$: \par
$(P1)$ take ${\hat \gamma }: {\hat M}\to N$ such that ${\hat \gamma
}=\gamma \circ \Xi $, then ${\bf P}_{{\hat \gamma },u}({\hat
s}_{0,q})=u_q$ for each $q=1,...,k$ and $\pi \circ {\bf P}_{{\hat
\gamma },u}={\hat \gamma }$
\par $(P2)$ ${\bf P}_{{\hat \gamma },u}$ is the $H^t_p$-mapping
by $\gamma $ and $u$;
\par $(P3)$ for each $x\in \hat M$ and every $\phi \in DifH^t_p
({\hat M}, \{ {\hat s}_{0,1},...,{\hat s}_{0,2k} \} )$ there is the
equality ${\bf P}_{{\hat \gamma },u}(\phi (x)) = {\bf P}_{{\hat
\gamma }\circ \phi ,u}(x)$, where $DifH^t_p({\hat M}, \{ {\hat
s}_{0,1},...,{\hat s}_{0,2k} \} )$ denotes the group of all $H^t_p$
homeomorphisms of $\hat M$ preserving marked points $\phi ({\hat
s}_{0,q})={\hat s}_{0,q}$ for each $q=1,...,2k$;
\par $(P4)$ ${\bf P}_{{\hat \gamma },u}$ is $G$-equivariant,
which means that ${\bf P}_{{\hat \gamma },uz}(x) = {\bf P}_{{\hat
\gamma },u}(x)z$ for every $x\in {\hat M}$ and each $z\in G$;
\par $(P5)$ if $U$ is an open neighborhood of ${\hat s}_{0,q}$ in
$\hat M$ and ${\hat \gamma }_0, {\hat \gamma }_1: U\to N$ are
$H^{t'}_p$-mappings such that ${\hat \gamma }_0({\hat
s}_{0,q})={\hat \gamma }_1({\hat s}_{0,q})=v_q$ and tangent spaces,
which are vector manifolds over ${\cal A}_r$, for $\gamma _0$ and
$\gamma _1$ at $v_q$ are the same, then the tangent spaces of ${\bf
P}_{{\hat \gamma }_0,u}$ and ${\bf P}_{{\hat \gamma }_1,u}$ at $u_q$
are the same, where $q=1,...,k$, $u=(u_1,...,u_k)$.

\par Two $H^{t'}_p$-differentiable principal $G$-bundles
$E_1$ and $E_2$ with parallel transport structures $(E_1,{\bf P}_1)$
and $(E_2,{\bf P}_2)$ are called isomorphic, if there exists an
isomorphism $h: E_1\to E_2$ such that ${\bf P}_{2,{\hat \gamma
},u}(x) = h ({\bf P}_{1,{\hat \gamma }, h^{-1}(u)}(x))$ for each
$H^t_p$-mapping $\gamma : M\to N$ and $u_q\in (E_2)_{y_0}$, where
$q=1,...,k$, $h^{-1}(u)=(h^{-1}(u_1),...,h^{-1}(u_k))$.

\par Let $(S^ME)_{t,H}:=
(S^{M, \{ s_{0,q}: q=1,...,k \} }E; N,G,{\bf P})_{t,H}$ be a set of
$H^t_p$-closures of isomorphism classes of $H^t_p$ principal $G$
fiber bundles with parallel transport structure.

\section{Structure of wrap groups}

\par {\bf 1. Proposition.} {\it The $H^m_p$ uniformity in
$L(S^m,N)$ (see \S 2.10 in \cite{luwrgfbqo}) for $m>1$ is strictly
stronger, than the $m$ times iterated $H^1_p$ uniformity.}
\par {\bf Proof.} If $f\in H^m$, then $\partial ^k f(x)/\partial
x_1^{k_1}...\partial x_m^{k_m}\in L^2$ for each $0\le k\le m$,
$k=k_1+...+k_m$, $0\le k_j$, $j=1,...,m$. But $g$ of $m$ times
iterated $H^1$ uniformity means that $\partial ^k g(x)/\partial
x_1^{k_1}...\partial x_m^{k_m}\in L^2$ for each $0\le k\le m$,
$k=k_1+...+k_m$, $0\le k_j\le 1$, $j=1,...,m$. The latter conditions
are weaker than that of $H^m$. For $m>1$ there may appear $g$ for
which such partial derivatives are not in $L^2$, when $1<k_j\le m$.
Using transition mappings of charts of atlases $At (M)$ and $At (N)$
and applying this locally we get the statement.

\par {\bf 2. Theorem.} {\it For a wrap group $W= (W^ME)_{t,H}$
(see Definition 2.7 \cite{luwrgfbqo}) there exists a skew product
${\hat W} = W{\tilde \otimes } W$ which is an $H^l_p$ alternative
Lie group and there exists a group embedding of $W$ into ${\hat W}$,
where $l=t'-t$ ($l=\infty $ for $t'=\infty $), $E=E(N,G,\pi ,\Psi )$
is a principal $G$-bundle of class $H^{t'}_p$ with $t'\ge t\ge [dim
(M)/2]+1$. If $G$ is associative, then $\hat W$ is associative.
Moreover, the loop group $L(S^1,E)$ is $H^t_p$ isomorphic with
$({\hat W}^{S^1}E)_{t,H}$ in the particular case of $S^1$.}

\par {\bf Proof.} Let ${\tilde W}$  be a set of all elements
$(g_1a_1\otimes g_2a_2)\in (W\otimes B)^2$, where $B$ is a free
non-commutative associative group with two generators $a, b$, $ab\ne
ba$, $g_1, g_2\in W$. Take in $\tilde W$ the equivalence relation:
$g_1g_2a\otimes g_2b {\tilde =} ~ g_1e_B\otimes ee_B$, for each
$g_1, g_2\in W$, where $e$ and $e_B$ denote the unit elements in $W$
and in $B$.

Define in $\tilde W$ the multiplication:
\par $(g_1a_1\otimes g_2a_2){\tilde \otimes }
(g_3a_3\otimes g_4a_4) := ((g_1g_3)(a_1a_3)
\otimes (g_4g_2)((a_1^{-1}a_4a_1) a_2)$ \\
for each $g_1, g_2, g_3, g_4\in W$ and every $a_1, a_2, a_3, a_4\in
B$, hence \par $(e\otimes g_1a_1){\tilde \otimes } (e\otimes g_2a_2)
= e\otimes (g_2g_1)(a_2a_1)$, \par $(g_1a_1\otimes e){\tilde \otimes
} (g_2a_2\otimes e) = (g_1g_2)(a_1a_2)\otimes e$,\par
$(g_1a_1\otimes e){\tilde \otimes } (e\otimes g_4a_4) =
g_1a_1\otimes g_4(a_1^{-1}a_4a_1)$, \par $(e\otimes g_4a_4){\tilde
\otimes } (g_1a_1\otimes e) := g_1a_1\otimes g_4a_4$. \\ Thus this
semidirect product $\tilde W$ of groups $(W\otimes B)\otimes ^s
(W\otimes B)$ is non-commutative, since $b^{-1}aba^{-1}\ne e$, where
$e := e\times e_B$, $\otimes ^s$ denotes the semidirect product,
$\otimes $ denotes the direct product.

Consider the minimal closed subgroup $A$ in the semidirect product
$\tilde W$ generated by elements $(g_1g_2a\otimes g_2b) {\tilde
\otimes } (g_1e_B\otimes ee_B)^{-1}$, where $B$ is supplied with the
discrete topology and $\tilde W$ is supplied with the product
uniformity. Then put ${\hat W} := {\tilde W}/A =: W{\tilde \otimes
}W$ and denote the multiplication in $\hat W$ as in $\tilde W$.

Therefore, $W$ has the group embedding $\theta : g\mapsto
(ge_B\otimes e)$ into ${\hat W}$ and the multiplication $m
[(g_1e_B\otimes e), (g_2e_B\otimes e)] = (g_1e_B\otimes e){\tilde
\otimes } (g_2e_B\otimes e)]$.

On the other hand, $(ga_1\otimes e) {\tilde \otimes } (e\otimes
ga_1a_2a_1^{-1}) = ga_1\otimes ga_2 = (e\otimes e) =: {\tilde e}$,
${\hat e} = {\tilde e}A=A$ is the unit element in $\hat W$ and
$(e\otimes ga_1a_2a_1^{-1}) = (ga_1\otimes e)^{-1}$ is the inverse
element of $(ga_1\otimes e)$, where $a_2\in B$ is such that
$(a_1\otimes a_2){\tilde \otimes }A = (e\otimes e){\tilde \otimes
}A=A$ in $\hat W$, $a_1=ea_1$, that is $a_1\otimes a_2 {\tilde =}
e\otimes e$ in $\tilde W$.

From preceding formulas it follows, that $\hat W$ is noncommutative
and alternative. As the manifold $\hat W$ is the quotient of the
$H^t_p$ manifold $W^2$ by the $H^t_p$ equivalence relation, hence
$\hat W$ is the $H^t_p$ differentiable space, since Conditions
$(D1-D4)$ of \S 2.1.3.2 \cite{luwrgfbqo} are satisfied. The group
operation and the inversion in $\hat W$ combines the product in $W$
and the inversion with the tensor product and the equivalence
relation, hence they are $H^l_p$ differentiable with $l=t'-t$,
$l=\infty $ for $t'=\infty $, (see \S \S 1.11, 1.12, 1.15 in
\cite{souriau} and \S 2.1.3.1 in \cite{luwrgfbqo}).
\par Then $((g_1\otimes g_2){\tilde \otimes } (g_3\otimes g_4))
{\tilde \otimes }(g_5\otimes g_6) := ((g_1g_3)g_5\otimes
g_6(g_4g_2))$ and \par $(g_1\otimes g_2){\tilde \otimes }
((g_3\otimes g_4)){\tilde
\otimes } (g_5\otimes g_6)) := (g_1(g_3g_5)\otimes (g_6g_4)g_2)$. \\
Therefore, $\hat W$ is alternative, since $W$ is alternative (see
Theorem 2.6.1 \cite{luwrgfbqo}) and $B$ is associative. If $G$ is
associative, then $W$ is associative and $\hat W$ is associative.

\par Consider the commutator \par $[(g_1a_1\otimes g_2a_2){\tilde {\otimes }}
(g_3a_3\otimes g_4a_4)]{\tilde {\otimes }}[(g_1a_1\otimes
g_2a_2)^{-1} {\tilde {\otimes }}$ \\ $(g_3a_3\otimes g_4a_4)^{-1}] =
\{ ((g_1g_3)(a_1a_3)\otimes (g_4g_2)((a_1^{-1}a_4a_1)a_2)){\tilde
{\otimes }}$ \\ $ [(g_1^{-1}a_1^{-1}\otimes
g_2^{-1}(a_1a_2^{-1}a_1^{-1})) {\tilde {\otimes
}}(g_3^{-1}a_3^{-1}\otimes g_4^{-1}(a_3a_4^{-1}a_3^{-1}))]$
\\ $=((g_1g_3)(a_1a_3)\otimes (g_4g_2)((a_1^{-1}a_4a_1)a_2)
{\tilde {\otimes }} ((g_1^{-1}g_3^{-1})(a_1^{-1}a_3^{-1})\otimes
(g_4^{-1}g_2^{-1})$
\\ $(a_1(a_3a_4^{-1}a_3^{-1})a_1^{-1})(a_1a_2^{-1}a_1^{-1})))
=(((g_1g_3)(g_1^{-1}g_3^{-1}))(a_1a_3a_1^{-1}a_3^{-1})\otimes $ \\
$((g_4^{-1}g_2^{-1})(g_4g_2)) ((a_1a_3)^{-1}
[((a_1a_3)a_4^{-1}(a_1a_3)^{-1})(a_1a_2^{-1}a_1^{-1})](a_1a_3))
((a_1^{-1}a_4a_1)a_2)$.

\par The minimal closed subgroup generated by products of such elements
is the commutant ${\tilde W}_c$ of $\tilde W$. The group
$(W^MN)_{t,H}$ is commutative (see Theorem 6$(2)$ \cite{luwrgfbqo}).
We have $B/B_c = \{ e \} $, the quotient group $G/G_c = G_{ab}$ is
the abelianization of $G$, particularly if $G$ is commutative, then
$G_{ab}=G$, where $G_c$ denotes the commutant subgroup of $G$.
Therefore, \par $(W^ME;N,G,{\bf P})_{t,H}/ [(W^ME;N,G,{\bf
P})_{t,H}]_c = (W^ME;N,G_{ab},{\bf P})_{t,H}$ \\ and inevitably
${\tilde W}/{\tilde W}_c = (W^ME;N,G_{ab},{\bf P})_{t,H}$.

\par Using the equivalence relation in ${\tilde W}$
we get ${\hat W}/{\hat W}_c = (W^ME;N,G_{ab},{\bf P})_{t,H}$.

\par In the particular case of $M=S^1$ for $g\in W$ take
$f\in g$, that is $<f>_{t,H} = g$. The equivalence class of $f$
relative to the analogous closures of orbits of the right action of
the subgroup $Diff^{\infty }_+(S^1,s_0)$ preserving a marked point
and an orientation of $S^1$ induced by that of $I=[0,1]$ denote by
$[f]_{t,H}$, then to $[f]_{t,H}$ put into the correspondence
$ga\otimes e$ in $\tilde W$, while to $[f^-]_{t,H}$ counterpose
$e\otimes gaba^{-1}$, where $f^-(x) := f(1-x)$ for each $x\in
[0,1]$, the unit circle $S^1$ is parametrized as $z = e^{2\pi i x}$,
$z\in S^1\subset \bf C$, $x\in [0,1]$. Their equivalence classes
$(ga\otimes e){\tilde \otimes }A$ and $(e\otimes gaba^{-1}){\tilde
\otimes }A$ in $\tilde W$ give elements in $\hat W$.

Since $[f]_{t,H}^{-1} := [f^-]_{t,H}$ and $[f_1\vee f_2]_{t,H} =
[f_1]_{t,H} [f_2]_{t,H}$, then $\hat W$ is isomorphic with
$L(S^1,E)_{t,H}$.

\par {\bf 3. Proposition.} {\it If there exists an
$H^{t'}_p$-diffeomorphism $\eta : N\to N$ such that $\eta (y_0) =
{y_0}'$, where $t\le t'$ then wrap groups $(W^ME; y_0)_{t,H}$ and
$(W^ME; {y_0}')_{t,H}$ defined with marked points $y_0$ and ${y_0}'$
are $H^l_p$-isomorphic as $H^l_p$-differentiable groups, where
$l=t'-t$ for finite $t'$, $l=\infty $ for $t'=\infty $.}
\par {\bf Proof.} Let $f\in H^t_p(M,E)$, then $\eta \circ \pi
\circ f(s_{0,q})= \eta (y_0)={y_0}'$ for each marked point $s_{0,q}$
in $M$, where $\pi : E\to N$ is the projection, $\pi \circ f=\gamma
$, $\gamma $ is a wrap, that is an $H^t_p$-mapping from $M$ into $N$
with $\gamma (s_{0,q})=y_0$ for $q=1,....,k$. The manifold $N$ is
connected together with $E$ and $G$ in accordance with conditions
imposed in \cite{luwrgfbqo}. Consider the $H^{t'}_p$-diffeomorphism
$\eta \times e$ of the principal bundle $E$. Then $\Theta :
H^t_p(M,W)\to H^t_p(M,W)$ is the induced isomorphism such that $\pi
\circ \Theta (f) := \eta \circ \pi \circ f: M\to N$ and $(\eta
\times e)\circ f = \Theta (f)$ for $f\in H^t_p(M,E)$. The mapping
$\Theta $ is $H^l_p$ differentiable by $f$, hence it gives the
$H^l_p$ isomorphism of the considered $H^l_p$-differentiable wrap
groups (see Theorem 6$(1)$ \cite{luwrgfbqo}).

\par {\bf 4. Remark.} As usually we suppose, that the principal bundle $E$,
its structure group $G$ and the base manifold $N$ are arcwise
connected. Let $({\cal P}^ME)_{t,H}$ be a space of equivalence
classes $<f>_{t,H}$ of $f\in H^t_p(M,W)$ relative to the closures of
orbits of the left action of $DifH^t_p(M; \{ s_{0,q}: q=1,...,k \}
)$. This means, that $({\cal P}^ME)_{t,H}$ is the quotient space of
$H^t_p(M,W)$ relative to the equivalence relation $R_{t,H}$.
\par There is the embedding $\theta : H^t_p(M,\{ s_{0,q}: q=1,...,k \} ; W)
\hookrightarrow H^t_p(M;W)$ and the evaluation mapping ${\hat ev}:
H^t_p(M;W)\to N^k$ such that ${\hat ev} (f) := ({\hat f}({\hat
s}_{0,q}): q=k+1,..,2k)$, ${\hat ev}_{{\hat s}_{0,q}} (f) := {\hat
f} ({\hat s}_{0,q})$, where ${\hat f} \in H^t_p({\hat M};W)$ is such
that ${\hat f} = f\circ \Xi $, $\Xi : {\hat M} \to M$ is the
quotient mapping. We get the diagram $H^t_p(M,\{ s_{0,q}: q=1,...,k
\} ; W)\to H^t_p(M;W)\to N^k$ with $H^t_p$ differentiable mappings,
which induces the diagram $H^{t,l+1}_p(M,\{ s_{0,q}: q=1,...,k \} ;
W, y_0)\to H^t_p(M,H^{t,l}_p(M,\{ s_{0,q}: q=1,...,k \} ;W, y_0)\to
H^{t,l}_p(M,\{ s_{0,q}: q=1,...,k \} ; W, y_0)$ for each $l\in \bf
N$, where $H^{t,l+1}_p(M,\{ s_{0,q}: q=1,...,k \} ; W, y_0) :=
H^t_p(M, \{ s_{0,q}: q=1,...,k \} ; H^{t,l}_p(M, \{ s_{0,q}:
q=1,...,k \} ; W, y_0))$, $H^{t,1}_p(M,\{ s_{0,q}: q=1,...,k \} ; W,
y_0) := H^t_p(M,\{ s_{0,q}: q=1,...,k \} ; W, y_0)$. Therefore,
there exist iterated wrap semigroups and groups $(S^ME)_{l+1;t,H} :=
(S^M(S^ME)_{l;t,H})_{t,H}$ and $(W^ME)_{l+1;t,H} :=
(W^M(W^ME)_{l;t,H})_{t,H}$, where $(S^ME)_{1;t,H}:= (S^ME)_{t,H}$
and $(W^ME)_{1;t,H} := (W^ME)_{t,H}$.
\par Evidently, if there are $H^t_p$ and $H^{t'}_p$
diffeomorphisms $\rho : M\to M_1$ and $\eta : N\to N_1$ mapping
marked points into respective marked points, then $H^t_p(M,W)$ is
isomorphic with $H^t_p(M_1,W_1)$ and hence $(W^ME)_{b;t,H}$ is
$H^t_p$ isomorphic as the $H^t_p$-manifold and $H^l_p$-isomorphic as
the $H^l_p$-Lie group with $(W^{M_1}E_1)_{b;t,H}$ for each $b\in \bf
N$, where $l=t'-t$, $l=\infty $ for $t'=\infty $, $t'\ge t\ge [dim
(M)/2]+1$. If $f: N\to N_1$ is a surjective map and $N$ is an
$H^t_p$-differentiable space, then $N$ inherits a structure of an
$H^t_p$-differentiable space with plots having the local form
$f\circ \rho : U\to N_1$, where $\rho : U\to N$ is a plot of $N$.

\par {\bf 5. Lemma.} {\it Let $E$ be an $H^{t'}_p$ principal bundle
and let $D$ be an everywhere dense subset in $N$ such that for each
$y\in D$ there exists an open neighborhood $V$ of $y$ in $N$ and a
differentiable map $p: V\to H^t_p(M, \{ s_{0,q}: q=1,...,k \} ;V,y)
:= \{ f\in H^t_p(M;V): f(s_{0,q}) =y, q=1,...,k \} $ such that
${\hat ev} _{{\hat s}_{0,q}}({\hat p}(y))=y$ for each $q=1,...,2k$
and each $y\in N$, where $p\circ \Xi = \pi \circ {\hat p}$. Then
${\hat ev}: H^t_p(M;W)\to N^k$ is an $H^t_p$ differentiable
principal $(S^ME)_{t,H}$ bundle.}

\par {\bf Proof.} Let $ \{ (V_j, y_j): j\in J \} $ be a family such that
$y_j\in V_j\cap D$ for each $j$ and there exists $p_j: V_j\to
H^t_p(M, \{ s_{0,q}: q=1,...,k \} ;V_j,y_j)$ so that ${\hat p}_j(
{\hat s}_{0,q})(y) =y\times e$ for each $q=1,...,2k$ and every $j$,
where $ \{ V_j: j\in J \} $ is an open covering of $N$, $y$ is a
constant mapping from $\hat M$ into $V_j$ with $y({\hat M})= \{ y \}
$, where ${\hat p}_j({\hat s}_{0,q})$ is the restriction to $V_j$ of
the projection ${\hat p}({\hat s}_{0,q}): ({\cal P}^ME)_{t,H}\to E$,
while $p_j(\Xi ({\hat x}))(y) = \pi \circ {\hat p}_j({\hat
x})(y\times e)$ for each $y\in N$ and $x= \Xi ({\hat x})$ in $M$,
where ${\hat x} \in \hat M$, $\Xi : {\hat M}\to M$. Then
$(W^ME)_{t,H}$ and $({\cal P}^ME)_{t,H}$ are supplied with the
$H^t_p$-differentiable spaces structure (see Remark 4 above and
Theorem 6 \cite{luwrgfbqo}), where the embedding
$(S^ME)_{t,H}\hookrightarrow ({\cal P}^ME)_{t,H}$ and the projection
${\hat ev}_{{\hat s}_{0,q}}: ({\cal P}^ME)_{t,H}\to N$ are
$H^t_p$-maps.
\par Let $\psi _j\in DifH^t_p(N)$ such that $\psi _j(y)=y_j$.
Specify a trivialization $\phi _j: {\hat p}_j^{-1}( {\hat s}_{0,q})
(V_j)\to V_j\times (S^ME)_{t,H}$ of the restriction ${\hat p}_j(
{\hat s}_{0,q})|_{V_j}$ of the projection ${\hat p}_j( {\hat
s}_{0,q}): ({\cal P}^ME)_{t,H}\to E$ by the formula $\phi _j(f) =
(f({\hat s}_{0,q}), \psi _j\circ {\hat p}_j({\hat s}_{0,q})(f))$ for
each $f \in ({\cal P}^ME)_{t,H}$ with $\pi \circ f ({\hat
s}_{0,q})=y$, where $\psi _j\circ {\hat p}_j(f) = \psi _j({\hat
p}_j(f))$. Then $\phi _j^{-1}(y,g) = g^{-1}(\psi _j\circ {\hat
p}_j(y))=:\eta $, $\eta \in ({\cal P}^ME)_{t,H}$ with $\pi \circ
\psi_j \circ f({\hat s}_{0,q})=y_j$, since $G$ is a group, where $g=
\psi _j\circ {\hat p}_j(f)$. Finally the combination of the family
$\{ {\hat ev} _{{\hat s}_{0,q}}: q=k+1,...,2k \} $ induce the
mapping ${\hat ev}: H^t_p(M;W)\to N^k$. By the construction a fiber
of this bundle is the monoid $(S^ME)_{t,H}$.

\par {\bf 6. Theorem.} {\it If $N$ is a smooth manifold
over ${\cal A}_r$ (holomorphic for $1\le r\le 3$ respectively), then
there exists an $H^t_p$-differentiable principal $(S^ME)_{t,H}$
bundle ${\hat ev}: ({\cal P}^ME)_{t,H}\to N^k$.}
\par {\bf Proof.} In view of Lemma 5 it is sufficient to prove that
for each $y\in N$ there exists a neighborhood $U$ of $y$ in $N$ and
an $H^t_p$-map $p_q: U\to H^t_p(M,W)$ such that
$ev_{s_{0,q}}(p_q(z))=z$ for each $q=1,...,k$, $z\in U$, where
$ev_x(f) = f(x)$.
\par  In $\hat M$ consider a rectifiable curve
$\zeta _q: [0,1]\to \hat M$ joining ${\hat s}_{0,q}$ with ${\hat
s}_{0,q+k}$, where $1\le q\le k$. Then consider a coordinate system
$(x_1,...,x_m)$ in $\hat M$ such that $x_1$ corresponds to a natural
coordinate along $\zeta _q$. This coordinate system is defined
locally for each chart of $\hat M$ and $x_1$ is defined globally.
\par Consider a real shadow $N_{\bf R}$ of $N$, then $N_{\bf R}$
is the Riemann $C^{\infty }$ manifold. Thus there exists a
Riemannian metric $\sf g$ in $N$. For each $y\in N$ there exists a
geodesic ball $U$ at $y$ of radius less than the injectivity radius
$\exp ^N$ for $\sf g$. Then there exists a map $p_q: U\to ({\cal
P}^MU)_{t,H}$ with $\pi \circ [p_q({\hat s}_{0,q+k})(z)] = z$ and
$\pi \circ [p_q({\hat s}_{0,q}(z)]=y$ for each $z\in U$, where
$p_q\circ \zeta _q=: {\hat \gamma }_{q,y,z}$ is the shortest
geodesic in $U$ joining $y$ with $z$, ${\hat \gamma }_{q,y,z}:
[0,1]\to N$, ${\hat \gamma }_{q,y,z}\circ \zeta ^{-1}_q(x_1)\in N$
for each $x_1$. Having initially ${\hat \gamma }_{q,y,z}$ extend it
to ${\hat p}_q$ on $\hat M$ with values in $E$ such that $p_q\circ
\Xi = \pi \circ {\hat p}_q$.

\par {\bf 7. Proposition.} {\it $(1)$.
The wrap group $(W^ME;N,G,{\bf P})_{t,H}$ is the principal $G^k$
bundle over $(W^MN)_{t,H}$. \par $(2)$. The abelianization
$[(W^ME;N,G,{\bf P})_{t,H}]_{ab}$ of the wrap group $(W^ME;N,G,{\bf
P})_{t,H}$ is isomorphic with $(W^ME;N,G_{ab},{\bf P})_{t,H}$. \par
$(3)$. For $n\ge 2$ the iterated loop group $(L^{S^n}E)_{t,H}$ is
isomorphic with the wrap group $(W^{S^n}E)_{t,H}$ for the sphere
$S^n$ and a principal fiber bundle $E$ for $dim_{\bf R}N\ge 2$ with
$k=1$.}

\par {\bf Proof. 1.} The bundle structure $\pi : E\to N$ induces
the bundle structure ${\hat {\pi }}: (W^ME;N,G,{\bf P})_{t,H}\to
(W^MN)_{t,H}$, since $\pi \circ {\bf P}_{{\hat \gamma }, u} ={\hat
\gamma }$. In view of Lemma 5 it is sufficient to show, that there
exists a neighborhood $U_G$ of $e$ in $(W^ME)_{t,H}$ and a
$G$-equivariant mapping $\phi : U_G\to (W^MN)_{t,H}$. Let $<{\bf
P}_{{\hat {\gamma }},u}>_{t,H}\in (W^ME)_{t,H}$, where ${\hat
{\gamma }}: {\hat M}\to N$, ${\hat {\gamma }} = \gamma \circ \Xi $,
$\gamma : M\to N$, $\gamma (s_{0,q})=y_0$ for each $q=1,...,k$. Then
$\pi \circ {\bf P}_{{\hat {\gamma }},u}={\hat {\gamma }}$ and ${\bf
P}_{{\hat {\gamma }},u}$ is $G$-equivariant by the conditions
defining the parallel transport structure, that is ${\bf P}_{{\hat
{\gamma }},u}(x)z= {\bf P}_{{\hat {\gamma }},uz}(x)$ for each $x\in
\hat M$ and $z\in G$ and every $u\in E_{y_0}$. We have that $uG =
\pi ^{-1}(y)$ for each $u\in E_y$ and $y\in N$.
\par Therefore, put $\phi = \pi _*$, where
$\pi _*<{\bf P}_{{\hat {\gamma }},u}>_{t,H} =<{{\hat {\gamma
}},u}>_{t,H}$ and take $U_G = \pi _*^{-1}(U)$, where $U$ is a
symmetric $U^{-1}=U$ neighborhood of $e$ in $(W^MN)_{t,H}$.

\par The group $G$ acts effectively on $E$. Since $G$ is arcwise
connected, then $G^k$ acts effectively on $(W^ME)_{t,H}$. Indeed,
for each $\zeta _q$ from \S 6 there is $g_q\in G$ corresponding to
${\hat \gamma }({\hat s}_{0,q+k})$ with ${\bf P}_{{\hat p}_q,{\hat
s}_{0,q}\times e}({\hat s}_{0,q+k}) =\{ y_0\times g_q \} \in
E_{y_0}$, $g_q\in G$ for every $q=1,...,k$. Moreover, $\pi
_*^{-1}(\pi _*(<{\bf P}_{{\hat {\gamma }},u}>_{t,H})) = <{\bf
P}_{{\hat {\gamma }},u}>_{t,H}G^k $. Then the fibre of ${\hat \pi }:
(W^ME;N,G,{\bf P})_{t,H}\to (W^MN)_{t,H}$ is $G^k$. Due to
Conditions 2$(P1-P5)$ \cite{luwrgfbqo} it is the principal $G^k$
differentiable bundle of class $H^t_p$.

\par {\bf 2, 3.} In view of Proposition 1 the loop
group $(L^{S^n}E)_{l,H}$ is everywhere dense in the $n$ times
iterated loop group $(L^{S^1}(...(L^{S^1}E)_{1,H}...)_{1,H}$, while
the wrap group $(W^{S^n}E)_{l,H}$  is everywhere dense in the $n$
times iterated wrap group $(W^{S^1}E)_{n;1,H}$ for each $l\ge n$.
For each $n>m$ there exists the natural projection $\pi ^m_n: S^n\to
S^m$ which induces the embeddings $(W^{S^m}E)_{t,H}\hookrightarrow
(W^{S^n}E)_{t,H}$ and $(L^{S^m}E)_{t,H}\hookrightarrow
(L^{S^n}E)_{t,H}$ in accordance with Corollary 9 \cite{luwrgfbqo},
since $k=1$ and choosing a marked point $s_0\in S^1$. Therefore, due
to $dim_{\bf R}N\ge 2$ the considered here wrap and loop groups are
infinite dimensional. Therefore, statements $(2,3)$ follow from
$(1)$ and the proof of Theorem 2 above and Proposition 11
\cite{luwrgfbqo} in accordance with which the iterated loop group
$(L^{S^1}(...(L^{S^1}E)_{1,H}...)_{1,H}$ is commutative.

\par {\bf 8. Proposition.} {\it If $E$ is contractible, then
$({\cal P}^ME)_{t,H}$ is contractible.}
\par {\bf Proof.} Let $g: [0,1]\times E\to E$ be a contraction
such that $g$ is continuous and $g(0,z)=z$ and $g(1,z)=y_0\times e$
for each $z\in E$. Then for each $f\in H^t_p(M,W)$ we get $g(0,f(x))
= f(x)$ and $g(1,f(x)) = y_0\times e$ for each $x\in M$. Moreover,
$g(s,<f>_{t,H}) \subset <g(s,f)>_{t,H}$ for each $s\in [0,1]$, since
$f\in g^{-1}_s(<g(s,f)>_{t,H})$ and $g$ is continuous while
$<g(s,f)>_{t,H}$ by its definition is closed in $H^t_p(M,W)$, where
$g_s(z) := g(s,z)$. Therefore, $id = g(0,*): ({\cal P}^ME)_{t,H}\to
({\cal P}^ME)_{t,H}$ and $g(1,({\cal P}^ME)_{t,H}) = <w_0>_{t,H}$.
\par {\bf 8.1. Notation.} Denote by $Hom^t_p((W^ME)_{t,H},G)$ or
$Hom^t_p((S^ME)_{t,H},G)$ the group or the monoid of $H^t_p$
differentiable homomorphisms from $(W^ME)_{t,H}$ or $(S^ME)_{t,H}$
respectively into $G$. By ${\cal A}_r^*$ is denoted the
multiplicative group of ${\cal A}_r\setminus \{ 0 \} $, where $0\le
r\le 3$.

\par {\bf 9. Theorem.} {\it Let $DifH^{t'}_p(N)$ acts transitively on $N$,
$t\le t'$. For each $H^{\infty }$ manifold $N$ and an $H^t_p$
differentiable group $G$ such that ${\cal A}_r^*\subset G$ with
$1\le r\le 3$ there exists a homomorphism of the $H^t_p$
differentiable space of all equivalence classes of $({\cal
P}^ME)_{t,H}$ relative to $DifH^{t'}_p(N)$ (see \S \S 1.3.2 and 3
\cite{luwrgfbqo}) and $Hom^t_p((S^ME)_{t,H},G^k)$. They are
isomorphic, when $G$ is commutative.}

\par {\bf Proof.} Mention that due to Theorem 6 the
$H^t_p$-differentiable principal $(S^ME)_{t,H}$ bundle ${\hat ev}:
({\cal P}^ME)_{t,H}\to N^k$ has a parallel transport structure
${\hat {\bf P}}_{{\hat \gamma },uz} (x)={\hat {\bf P}}_{{\hat
{\gamma }},u}(x)z$ for each $x\in \hat M$ and all $\gamma \in
H^t_p(M,N)$ and $u\in {\hat ev}^{-1}(\gamma (s_{0,k}))$ and every
$z\in G$ and the corresponding ${\hat \gamma }: {\hat M}\to N$ such
that $\gamma \circ \Xi = \hat {\gamma }$. If $x = {\hat s}_{0,q}$
with $1\le q\le k$, then ${\hat {\bf P}}$ gives the identity
homomorphism from $(S^ME)_{t,H}$ into $(S^ME)_{t,H}$. If $\theta :
(S^ME)_{t,H}\to G^k$ is an $H^t_p$ differentiable homomorphism, then
the holonomy of the associated parallel transport ${\hat {\bf
P}}^{\theta }$ on the bundle $({\cal P}^ME)_{t,H}\times
^{\theta}G\to N^k$ is the homomorphism $\theta : (S^ME)_{t,H}\to
G^k$ (see \S 2.3 in \cite{luwrgfbqo}). At the same time the group
$G$ contains continuous one-parameter subgroups from ${\cal A}_r^*$,
where $1\le r\le 3$. If $g\in (W^MN)_{t,H}$ and $g\ne e$, then $g$
is of infinite order, since $w_0$ does not belong to $g^n$ for each
$n\ne 0$ non-zero integer $n$, where $w_0(M)= \{ y_0 \} $.

\par This holonomy induces a map $h: ({\cal P}^ME)_{t,H}/{\cal Q} \to
Hom^t_p((S^ME)_{t,H},G^k)$, where ${\cal Q}$ is an equivalence
relation caused by the transitive action of $DifH^{t'}_p(N)$ such
that $(S^ME)_{t,H}$ with distinct marked points either $ \{ s_{0,q}:
q=1,...,k \} $ in $M$ and $y_0$ or ${\tilde y}_0$ in $N$ are
isomorphic, since there exists $\psi \in DifH^{t'}_p(N)$ such that
$\psi (y_0)={\tilde y}_0$.
\par If $G$ is commutative, then this map is the homomorphism, since
$(S^ME)_{t,H}$ is the commutative monoid for a commutative group $G$
(see Theorem 3.2 \cite{luwrgfbqo}) and $u{\bf P}_{{\hat {\gamma
}}_1,v_1}(x_1) {\bf P}_{{\hat {\gamma }}_2,v_2}(x_2)= u{\bf
P}_{{\hat {\gamma }}_2,v_2}(x_2) {\bf P}_{{\hat {\gamma
}}_1,v_1}(x_1)$ for each $x_1, x_2\in \hat M$ and $u, v_1, v_2 \in
E_{y_0}$. There is the embedding $(S^ME)_{t,H}\hookrightarrow
(W^ME)_{t,H}$, hence a homomorphism $\theta : (W^ME)_{t,H}\to G^k$
has the restriction on $(S^ME)_{t,H}$ which is also the
homomorphism. \par For $G\supset {\cal A}_r^*$ there exists a family
of $f\in Hom^t_p((S^ME)_{t,H},G^k)$ separating elements of the wrap
monoid $(S^ME)_{t,H}$, hence there exists the embedding of
$(S^ME)_{t,H}$ into $Hom^t_p((S^ME)_{t,H},G^k)$. The bundle $({\cal
P}^ME)_{t,H}\times ^{\theta }G\to N^k$ has the induced parallel
transport structure ${\bf P}^{\theta }$. The holonomy of the
parallel transport structure on $({\cal P}^MN)_{t,H}\times ^{\theta
}G\to N^k$ is $\theta $. Therefore, the map
$H^t_p((S^ME)_{t,H},G^k)\ni \theta \mapsto {\bf P}^{\theta }$ is
inverse to $h$.

\par {\bf 10. Theorems.} {\it Suppose that $M_2\hookrightarrow M_1$
and $M = M_1\setminus (M_2\setminus \partial M_2)$ and ${\hat
M}_2\hookrightarrow {\hat M}_1$ and ${\hat M} = {\hat M}_1\setminus
({\hat M}_2\setminus \partial {\hat M}_2)$ and $N_2\hookrightarrow
N_1$ are $H^t_p$-pseudo-manifolds with the same marked points $\{
s_{0,q}: q=1,...,k \} $ for $M_1$ and $M_2$ and $M$ and $y_0\in N_2$
satisfying conditions of \S 2 \cite{luwrgfbqo} and $G_2$ is a closed
subgroup in $G_1$ with a topologically complete principal fiber
bundle $E$ with a structure group $G_1$.
\par {\bf 1}. Then $(W^{M_2, \{ s_{0,q}: q=1,...,k \} }E;N_2,G_2,{\bf
P})_{t,H}$ has an embedding as a closed subgroup into $(W^{M_1, \{
s_{0,q}: q=1,...,k \} } E;N_1,G_1,{\bf P})_{t,H}$.
\par {\bf 2}. The wrap group
$(W^{M_2, \{ s_{0,q}: q=1,...,k \} }E;N,G_2,{\bf P})_{t,H}$ is
normal in \\ $(W^{M_1, \{ s_{0,q}: q=1,...,k \} } E;N,G_1,{\bf
P})_{t,H}$ if and only if $G_2$ is a normal subgroup in $G_1$.
\par {\bf 3}. In the latter case $(W^ME;N,G,{\bf P})_{t,H}$
is isomorphic with \\ $(W^{M_1}E;N,G_1,{\bf
P})_{t,H}/(W^{M_2}E;N,G_2,{\bf P})_{t,H}$, where $G=G_1/G_2$.}

\par {\bf Proof. 1.} If ${\hat {\gamma }}_2\in H^t_p({\hat M}_2, N_2)$,
then it has an $H^t_p$ extension to ${\hat {\gamma }}_1\in
H^t_p({\hat M}_1, N_1)$ due to Theorem III.4.1 \cite{miha}.
Therefore, the parallel transport structure ${\bf P}_{{\hat {\gamma
}}_1,u}$ over ${\hat M}_1$ serves as an extension of ${\bf P}_{{\hat
{\gamma }}_2,u}$ over ${\hat M}_2$. The uniform spaces $H^t_p(M_j,
\{ s_{0,1},...,s_{0,k} \} ;W_j,y_0)$ are complete for $j=1,2$, since
the principal fiber bundle $E$ is topologically complete and the
corresponding principal fiber sub-bundle $E_2$ with the structure
group $G_2$ is also complete (see Theorem 8.3.6 \cite{eng}).
Therefore, $H^t_p(M_2, \{ s_{0,1},...,s_{0,k} \} ;W_2,y_0)$ has
embedding as the closed subspace into $H^t_p(M_1, \{
s_{0,1},...,s_{0,k} \} ;W_1,y_0)$. Each $H^t_p$ diffeomorphism of
$M_2$ has an $H^t_p$ extension to a diffeomorphism of $M_1$ (see
also \S III.4 in \cite{miha} and \cite{touger}). Since $G_2$ is a
closed subgroup in $G_1$, then $(S^{M_2, \{ s_{0,q}: q=1,...,k \}
}E;N_2,G_2,{\bf P})_{t,H}$ has an embedding as a closed sub-monoid
into $(S^{M_1, \{ s_{0,q}: q=1,...,k \} } E;N_1,G_1,{\bf P})_{t,H}$
and inevitably $(W^{M_2, \{ s_{0,q}: q=1,...,k \} }E;N_2,G_2,{\bf
P})_{t,H}$ has an embedding as a closed subgroup into $(W^{M_1, \{
s_{0,q}: q=1,...,k \} } E;N_1,G_1,{\bf P})_{t,H}$ due to Theorem 6.1
\cite{luwrgfbqo}.

\par {\bf 2.} The groups $(W^{M_j, \{ s_{0,q}: q=1,...,k \} }N)_{t,H}$
for $j=1, 2$ are commutative and $(W^{M_j, \{ s_{0,q}: q=1,...,k \}
}E)_{t,H}$ is the $G_j^k$ principal fiber bundle on $(W^{M_j, \{
s_{0,q}: q=1,...,k \} }N)_{t,H}$ (see Theorem 6.2 \cite{luwrgfbqo}
and Proposition 7.1 above). Therefore, $(W^{M_2, \{ s_{0,q}:
q=1,...,k \} }E)_{t,H}$ is the normal subgroup in $(W^{M_1, \{
s_{0,q}: q=1,...,k \} }E)_{t,H}$ if and only if $G_2$ is the normal
subgroup in $G_1$.

\par {\bf 3.} Consider the principal fiber bundle $E(N,G,\pi ,\Psi )$
with the structure group $G$ (see Note 1.3.2 \cite{luwrgfbqo}) and
the parallel transport structure ${\bf P}$ for the $H^t_p$
pseudo-manifold $\hat M$, where $G=G_1/G_2$ is the quotient group.
If ${\hat {\gamma }}_1\in H^t_p({\hat M}_1, N)$, then ${\hat {\gamma
}}_1$ is the combination
\par $(i)$ ${\hat {\gamma }}_1 = {\hat {\gamma }}_2\nabla {\hat
{\gamma }}$,
\\ where ${\hat {\gamma }}_2$ and ${\hat {\gamma }}$ are
restrictions of ${\hat {\gamma }}_1$ on ${\hat M}_2$ and ${\hat M}$
correspondingly. On the other hand, each ${\hat {\gamma }}\in
H^t_p({\hat M}, N)$ has an extension ${\hat {\gamma }}_1\in
H^t_p({\hat M}_1, N)$. The manifold ${\hat M}_1$ is metrizable by a
metric $\rho $. For each $\epsilon
>0$ there exists $\psi \in DifH^t_p({\hat M}_1; \{ {\hat s}_{0,q}:
q=1,...,2k \} )$ such that $(\psi ({\hat M})\cap {\hat M}_2)\subset
\bigcup_{l=1}^s B({\hat M}_1,x_l,\epsilon )$ for some $x_l\in {\hat
M}_1$ with $l=1,...,s$ and $s\in \bf N$ and $\psi |_{{\hat
M}_1\setminus ({\hat M}\bigcup_{l=1}^s B({\hat M}_1,x_l,\epsilon
))}=id$, since ${\hat M}_1$ and ${\hat M}_2$ are compact
pseudo-manifolds. Therefore, using Lemma 2.1.3.16 \cite{lufoclg} and
charts of the manifolds gives \par $<{\bf P}_{{\hat {\gamma
}},u}|_M>_{t,H} = <{\bf P}_{{\hat {\gamma
}}_1,u}|_{M_1}>_{t,H}/<{\bf P}_{{\hat {\gamma }}_2,u}|_{M_2}>_{t,H}$
\\ due to decomposition $(i)$, since ${\bf P}_{{\hat {\gamma
}},u}|_{M_j}\in G_j$ for $j=1,2$ and $G=G_1/G_2$ is the $H^{t'}_p$
quotient group with $t'\ge t$.
Consequently, $(W^ME;N,G,{\bf P})_{t,H}$ is isomorphic with \\
$(W^{M_1}E;N,G_1,{\bf P})_{t,H}/(W^{M_2}E;N,G_2,{\bf P})_{t,H}$ (see
also \S \S 3, 6 \cite{luwrgfbqo}).

\par {\bf 11. Corollary.} {\it Let suppositions of Theorem 10 be
satisfied. Then $(W^MN)_{t,H}$ is isomorphic with $(W^{M_1}N)_{t,H}/
(W^{M_2}N)_{t,H}$.}
\par {\bf Proof.} For $(W^MN)_{t,H}$ taking $G=G_1=G_2=\{ e \} $
we get the statement of this corollary from Theorem 10.3.

\par {\bf 12. Proposition.} {\it  Suppose that $M=M_1\vee M_2$, where
$M_1$ and $M_2$ are $H^t_p$-pseudo-manifolds satisfying Conditions
2.2$(i-v)$ \cite{luwrgfbqo} with the bunch taken by marked points
$\{ s_{0,q}: q=1,...,k \} $, then $(W^MN)_{t,H}$ is isomorphic with
the internal direct product $(W^{M_1}N)_{t,H}\otimes
(W^{M_2}N)_{t,H}$.}

\par {\bf Proof.} The manifold $M$ has marked points
$\{ s_{0,q}: q=1,...,k \} $ such that $s_{0,q}$ corresponds to
$s_{0,q,1}$ glued with $s_{0,q,2}$ in the bunch $M_1\vee M_2$ for
each $q=1,...,k$, where $s_{0,q,j}\in M_j$ are marked points $j=1,
2$. Since each $M_j$ satisfies Conditions 2.2$(i-v)$
\cite{luwrgfbqo}, then $M$ satisfies them also. \par In view of
Theorem 10.1 $(W^{M_j, \{ s_{0,q}: q=1,...,k \} }N)_{t,H}$ has an
embedding as a closed subgroup into $(W^{M, \{ s_{0,q}: q=1,...,k \}
} N)_{t,H}$ for $j=1, 2$. If $\gamma _j\in H^t_p(M_j, \{ s_{0,q}:
q=1,...,k \}; N,y_0)$ for $j=1, 2$, then $\gamma _1\vee \gamma _2\in
H^t_p(M, \{ s_{0,q}: q=1,...,k \}; N,y_0)$. On the other hand, each
$\gamma \in H^t_p(M, \{ s_{0,q}: q=1,...,k \}; N,y_0)$ has the
decomposition $\gamma = \gamma _1\vee \gamma _2$, where $\gamma _j =
\gamma |_{M_j}$ for $j=1, 2$. Therefore, $<\gamma >_{t,H} = <\gamma
_1\vee w_{0,2}>_{t,H}\vee <w_{0,1}\vee \gamma _2>_{t,H}$, where
$w_0(M)= \{ y_0 \} $, $w_{0,j}= w_0|_{M_j}$ for $j=1, 2$, hence
$(W^MN)_{t,H}$ is isomorphic with $(W^{M_1}N)_{t,H}\otimes
(W^{M_2}N)_{t,H}$.

\par {\bf 13. Propositions. 1.} {\it Let $\theta : N_1\to N$ be an
embedding with $\theta (y_1)=y_0$, or $F: E_1\to E$ be an embedding
of principal fiber bundles over ${\cal A}_r$ such that $\pi \circ
F|_{N_1\times e} =\theta \circ \pi _1$, then there exist embeddings
$\theta _*: (W^MN_1)_{t,H}\to (W^MN)_{t,H}$ and $F_*:
(W^ME_1)_{t,H}\to (W^ME)_{t,H}$.}
\par {\bf 2.} {\it If $\theta : N_1\to N$ and $F: E_1\to E$ are a quotient
mapping and a quotient homomorphism such that $N_1$ is a covering
pseudo-manifold of a pseudo-manifold $N$, then $(W^MN)_{t,H}$ is the
quotient group of some closed subgroup in $(W^MN_1)_{t,H}$ and
$(W^ME)_{t,H}$ is the quotient group of some closed subgroup in
$(W^ME_1)_{t,H}$.}
\par {\bf 3.} {\it If there are an $H^t_p$ diffeomorphism $f_1: M\to M_1$
and an $H^{t'}_p$-isomorphism $f_2: E\to E_1$, then wrap groups
$(W^{M_1}E_1)_{t,H}$ and $(W^ME)_{t,H}$ are isomorphic.}

\par {\bf Proof. 1.} If $\gamma _1\in H^t_p(M,\{ s_{0,q}: q=1,...,k \};
N_1,y_1)$, then $\theta \circ \gamma _1=\gamma \in H^t_p(M,\{
s_{0,q}: q=1,...,k \}; N,y_0)$, $<\gamma >_{t,H} = \theta _*<\gamma
_1>_{t,H}$, where $\theta _*<\gamma _1>_{t,H} := \{ \theta \circ f:
fR_{t,H}\gamma _1 \} $. In addition $F|_{E_{1,v}}$ gives an
embedding $F: G_1\to G$, where $G_1$ and $G$ are structural groups
of $E_1$ and $E$. Therefore, for the parallel transport structures
we get \par $(1)$ $F\circ {\bf P}^1_{{\hat {\gamma }}_1,v}(x) = {\bf
P}_{{\hat {\gamma }},u}(x)$\\ for each $x\in \hat M$, where
$F(v)=u$, $\pi \circ F = \theta \circ \pi _1$, where ${\bf P}^1$ is
for $E_1$ and $\bf P$ for $E$. Define $F_*<{\bf P}^1_{{\hat {\gamma
}}_1,v}>_{t,H} := \{ F\circ g: gR_{t,H}{\bf P}^1_{{\hat {\gamma
}}_1,v} \} $. Since $\theta $ and $F$ are $H^t_p$ differentiable
mappings, then $\theta _*$ and $F_*$ are embeddings of $H^t_p$
manifolds and group homomorphisms of $H^l_p$ differentiable groups
(see also Theorems 6 \cite{luwrgfbqo}).

\par {\bf 2.} If $\gamma \in H^t_p(M,\{ s_{0,q}: q=1,...,k \};
N,y_0)$, then there exists $\gamma _1 \in H^t_p(M,\{ s_{0,q}:
q=1,...,k \}; N_1,y_1)$ such that $\theta \circ \gamma _1=\gamma $,
since $N_1$ is a covering of $N$, that is each $y\in N$ has a
neighborhood $V_y$ for which $\theta ^{-1}(V_y)$ is a disjoint union
of open subsets in $N_1$ for each $y\in N$. This $\gamma _1$ exists
due to connectedness of $M$ and $\gamma (M)$, where $\gamma
(M)\subset N$. To each parallel transport in $E_1$ there corresponds
a parallel transport in $E$ so that Equation $(1)$ above is
satisfied. Put $\theta _*^{-1}<\gamma >_{t,H} = \{ <\gamma
_1>_{t,H}: \theta \circ \gamma _1 = \gamma \} $ and $F_*^{-1}<{\bf
P}_{{\hat {\gamma }},u}>_{t,H} := \{ <{\bf P}^1_{{\hat {\gamma
}}_1,v}>_{t,H}: F\circ {\bf P}^1_{{\hat {\gamma }}_1,v} = {\bf
P}_{{\hat {\gamma }},u} \} $, where $F(v)=u$.

This gives quotient mappings $\theta _*$ and $F_*$ from closed
subgroups $\theta _*^{-1}(W^MN)_{t,H}$ and $F_*^{-1}(W^ME)_{t,H}$ in
$(W^MN_1)_{t,H}$ and $(W^ME_1)_{t,H}$ respectively onto
$(W^MN)_{t,H}$ and $(W^ME)_{t,H}$ by closed subgroups $\theta
_*^{-1}(e)$ and $F_*^{-1}(e)$ correspondingly.

\par {\bf 3.} We have that $g \in H^t_p(M,\{ s_{0,q}: q=1,...,k \};W,y_0)$
if and only if $f_2\circ g \circ f_1^{-1}\in H^t_p(M_1,\{ s_{0,q,1}:
q=1,...,k;W_1,y_1)$, where $f_1(s_{0,q})=s_{0,q,1}$ for each
$q=1,...,k$, $f_2(y_0\times e)=y_1\times e$. At the same time $\psi
\in DifH^t_p(M)$ if and only if $f_1\circ \psi \circ f_1^{-1}\in
DifH^t_p(M_1)$. Hence $(S^ME)_{t,H}$ is isomorphic with
$(S^{M_1}E_1)_{t,H}$ and inevitably wrap groups $(W^ME)_{t,H}$ and
$(W^{M_1}E_1)_{t,H}$ are $H^t_p$ diffeomorphic as manifolds and
isomorphic as $H^l_p$ groups.

\par {\bf 14. Note.} If $N$ is a manifold not necessarily orientable,
then it contains up to equivalence of atlases a connected chart $V$
open in $N$ such that $y\in V$ and $V$ is orientable. Since
$(W^ME|_V)_{t,H}$ is the infinite dimensional group, then
$(W^ME)_{t,H}$ is also infinite dimensional even if $N$ is not
orientable due to Proposition 13.1. If $N$ is not orientable, then
there exists an orientable covering manifold $N_1$ and a quotient
mapping $\theta : N_1\to N$ as in Proposition 13$(2)$ (see also
about coverings and orientable coverings in \S \S 50, 51
\cite{pont}, \S \S II.4.18,19 \cite{dubnovfom}).
\par It is necessary to mention that some circumstances of wrap
groups are related also with their infinite dimensionality.

\par {\bf 15. Note.}
Let $G$ be a topological group not necessarily associative, but
alternative:
\par $(A1)$ $g(gf)=(gg)f$ and $(fg)g=f(gg)$ and $g^{-1}(gf)=f$
and $(fg)g^{-1}=f$ for each $f, g\in G$ \\
and having a conjugation operation which is a continuous
automorphism of $G$ such that \par $(C1)$ $conj (gf)=conj (f)
conj(g)$ for each $g, f\in G$, \par $(C2)$ $conj (e)=e$ for the unit
element $e$ in $G$.
\par If $G$ is of definite class of smoothness, for example, $H^t_p$
differentiable, then $conj$ is supposed to be of the same class. For
commutative group in particular it can be taken the identity mapping
as the conjugation. For $G= {\cal A}_r^*$ it can be taken $conj
(z)={\tilde z}$ the usual conjugation for each $z\in {\cal A}_r^*$,
where $1\le r\le 3$.
\par Suppose that \par $(A2)$ ${\hat G} =
{\hat G}_0i_0\oplus {\hat G}_1i_1\oplus ... \oplus {\hat
G}_{2^r-1}i_{2^r-1}$ such that $G$ is a multiplicative group of a
ring $\hat G$ with the multiplicative group structure, where ${\hat
G}_0,...,{\hat G}_{2^r-1}$ are pairwise isomorphic commutative
associative rings and $ \{ i_0,...,i_{2^r-1} \} $ are generators of
the Cayley-Dickson algebra ${\cal A}_r$, $1\le r\le 3$ and
$(y_li_l)(y_si_s)=(y_ly_s)(i_li_s)$ is the natural multiplication of
any pure states in $G$ for $y_l\in G_l$. For example, $G=({\cal
A}_r^*)^n$ and ${\hat G} = {\cal A}_r^n$.

\par {\bf 16. Lemma.} {\it If $G$ and $K$ are two topological
or differentiable groups twisted over $ \{ i_0,...,i_{2^r-1} \} $
satisfying conditions 15$(A1,A2,C1,C2)$ and $K$ is a closed normal
subgroup in $G$, where $2\le r\le 3$, then the quotient group is
topological or differentiable and twisted over $ \{
i_0,...,i_{2^r-1} \} $.}
\par {\bf Proof.} Since ${\hat G}= {\hat G}_0i_0\oplus {\hat G}_1i_1
\oplus ... \oplus {\hat G}_{2^r-1}i_{2^r-1}$, where ${\hat
G}_0,...,{\hat G}_{2^r-1}$ are pairwise isomorphic, then ${\hat
G}/{\hat K} = ({\hat G}_0/{\hat K}_0)i_0\oplus ... \oplus ({\hat
G}_{2^r-1}/{\hat K}_{2^r-1})i_{2^r-1}$ is also twisted. Each ${\hat
G}_j$ is associative, hence $G/K$ is alternative, since $2\le r\le
3$ and using multiplicative properties of generators of the
Cayley-Dickson algebra ${\cal A}_r$. On the other hand, $conj
(K)=K$, hence $conj (gK)=K conj (g)= conj (g) K\in G/K$ and $conj
(ghK)= conj (gh) K=(conj (h) conj (g)) K = (conj (h)K) (conj (g)K) =
conj (hK) conj (gK)= conj (gK hK)$. \par The subgroup $K$ is closed
in $G$, hence by the definition of the quotient differentiable
structure $G/K$ is the differentiable group (see also \S \S 1.11,
1.12, 1.15 in \cite{souriau}).

\par {\bf 17. Proposition.} {\it Let $\eta : N_1\to N_2$
be an $H^{t'}_p$-retraction of $H^{t'}_p$ manifolds, $N_2\subset
N_1$, $\eta |_{N_2}=id$, $y_0\in N_2$, where $t'\ge t$, $M$ is an
$H^t_p$ manifold, $E(N_1,G,\pi ,\Psi )$ and $E(N_2,G,\pi ,\Psi )$
are principal $H^{t'}_p$ bundles with a structure group $G$
satisfying conditions of \S 2 \cite{luwrgfbqo}. Then $\eta $ induces
the group homomorphism $\eta _*$ from $(W^ME;N_1,G,{\bf P})_{t,H}$
onto $(W^ME;N_2,G,{\bf P})_{t,H}$.}
\par {\bf Proof.} In view of Proposition 7$(1)$ the wrap group
$(W^ME;N_1,G,{\bf P})_{t,H}$ is the principal $G^k$ bundle over
$(W^MN_1)_{t,H}$. Extend $\eta $ to $\vartheta : E(N_1,G,\pi ,\Psi
)\to E(N_2,G,\pi ,\Psi )$ such that $\pi _2 \circ \vartheta =\eta
\circ \pi _1$ and $pr_2\circ \vartheta =id: G\to G$, where $pr_2:
E_y\to G$ is the projection, $y\in N_1$. If $f\in H^t_p(M,N_1)$,
then $\eta \circ f := \eta (f(*))\in H^t_p(M,N_2)$. If
$f(s_{0,q})=y_0$, then $\eta (f(s_{0,q}))=y_0$, since $y_0\in N_2$.
Since $N_2\subset N_1$, then $H^t_p(M,N_2)\subset H^t_p(M,N_1)$. The
parallel transport structure $\bf P$ is over the same manifold $M$.
\par Put $\eta _* (<{\bf P}_{{\hat {\gamma }},u}>_{t,H}) = <{\bf
P}_{\eta \circ {\hat {\gamma }},u}>_{t,H}$, where ${\hat {\gamma }}:
{\hat M}\to N_1$. In view of Theorems 2.3 and 2.6 \cite{luwrgfbqo}
$\eta _*(<{\bf P}_{{\hat {\gamma }}_1,u}\vee {\bf P}_{{\hat {\gamma
}}_2,u}>_{t,H} = \eta ^*(<{\bf P}_{{\hat {\gamma }}_1,u}>_{t,H})
\eta _*(<{\bf P}_{{\hat {\gamma }}_2,u}>_{t,H})$, and we can put
$\eta _*(q^{-1})= [\eta _*(q)]^{-1}$, consequently, $\eta _*$ is the
group homomorphism. Moreover, for each $g\in (W^ME;N_2,G,{\bf
P})_{t,H}$ there exists $q\in (W^ME;N_1,G,{\bf P})_{t,H}$ such that
$\eta _*(q)=g$, since $\gamma : M\to N_2$ and $N_2\subset N_1$ imply
$\gamma : M\to N_1$, while the structure group $G$ is the same,
hence $\eta _*$ is the epimorphism.

\par {\bf 18. Definition.} Let $G$ be a topological group
satisfying Conditions 15$(A1,A2,C1,C2)$ such that $G$ is a
multiplicative group of the ring $\hat G$, where $1\le r\le 2$. Then
define the smashed product $G^s$ such that it is a multiplicative
group of the ring ${\hat G}^s := {\hat G}\otimes _l{\hat G}$, where
$l=i_{2^r}$ denotes the doubling generator, the multiplication in
${\hat G}\otimes _l{\hat G}$ is \par $(1)$ $(a+bl)(c+vl) = (ac - v^*
b) + (va+bc^*)l$ for each $a, b, c, v\in \hat G$, where $v^* = conj
(v)$.
\par A smashed product $M_1\otimes _lM_2$ of manifolds $M_1, M_2$
over ${\cal A}_r$ with $dim (M_1)= dim (M_2)$ is defined to be an
${\cal A}_{r+1}$ manifold with local coordinates $z=(x,yl)$, where
$x$ in $M_1$ and $y$ in $M_2$ are local coordinates.

\par Its existence and detailed description are demonstrated below.

\par {\bf 19. Proposition.} {\it The ring ${\hat G}^s$
has a multiplicative group $G^s$ containing all $a+bl\ne 0$ with $a,
b\in  {\hat G}$. If $\hat G$ is a topological or $H^t_p$
differentiable ring over ${\cal A}_r$ for $t\ge dim (G)+1$, then
${\hat G}^s$ is a topological or $H^t_p$ differentiable over ${\cal
A}_{r+1}$ ring.}
\par {\bf Proof.} For each $1\le r \le 2$ the group $G$ is
associative, since the generators $ \{ i_0,...,i_{2^r-1} \} $ form
the associative group, when $r\le 2$. An element $a+bl\in {\hat
G}^s$ is non-zero if and only if $(a+bl)(a+bl)^* = aa^* +bb^*\ne 0$
due to 15$(A1,A2,C1,C2)$ and 18$(1)$. For $a+bl\ne 0$ put $u = (a^*
- lb^*)/(aa^*+bb^*)$, where $aa^*+bb^*\in G_0$, hence $u (a+bl) =
(a+bl) u= 1\in G_0$, since $G_j$ is commutative for each
$j=0,...,2^r-1$, where $G_j$ denotes the multiplicative group of the
ring ${\hat G}_j$. For $r\le 2$ the family of generators $ \{
i_0,..., i_{2^{r+1}-1} \} $ forms the alternative group, hence
${\hat G}^s = {\hat G}_0i_0\oplus ... \oplus {\hat G}_{2^{r+1}-1}$
is alternative, where ${\hat G}_j$ are isomorphic with ${\hat G}_0$
for each $j$.
\par If an addition in $\hat G$ is continuous, then evidently
$(a+bl) + (c+ql)= (a+c) +(b+q)l$ is continuous. If the
multiplication in $\hat G$ is continuous, then Formula 18$(1)$ shows
that the multiplication in ${\hat G}^s$ is continuous as well.
\par We have the decomposition ${\cal A}_{r+1} = {\cal A}_r\oplus
{\cal A}_rl$. If $\hat G$ is $H^t_p$ differentiable, then from the
definition of plots it follows, that ${\hat G}^s$ is $H^t_p$
differentiable over ${\cal A}_{r+1}$ (see also in details
20$(1-5)$).

\par {\bf 20. Theorem.} {\it Let $M_1, M_2$ and $N_1, N_2$ be
$H^t_p$ manifolds over ${\cal A}_r$ with $1\le r \le 2$, and let $G$
be a group satisfying Conditions 15$(A1,A2,C1,C2)$, let also
$M_1\otimes _lM_2$, $N_1\otimes _lN_2$ be smashed products of
manifolds and $G^s$ be a smashed product group (see Proposition 19),
where $dim (M_1) = dim (M_2)$, $dim (N_1) = dim (N_2)$, $t\ge \max
(dim (M_1), dim (N_1), dim (G)) +1$. Then the wrap group
$(W^{M_1\otimes _lM_2; \{ s_{0,j,1}\otimes _l s_{0,v,2}:
j=1,...,k_1; v=1,...,k_2 \} } E;N_1\otimes _lN_2,G^s,{\bf
P}^s)_{t,H}$ is twisted over $ \{ i_0,...,i_{2^{r+1}-1} \} $ and is
isomorphic with the smashed product
\par $W^{M_2; \{ s_{0,v,2}: v=1,...,k_2 \} }E;N_1,(W^{M_1; \{
s_{0,j,1}: j=1,...,k_1 \} }E;N_1,G,{\bf P}_1)_{t,H},{\bf
P}_2)_{t,H}\otimes _l$
\par $W^{M_2; \{ s_{0,v,2}: v=1,...,k_2 \} }E;N_2,(W^{M_1; \{
s_{0,j,1}: j=1,...,k_1 \} }E;N_2,G,{\bf P}_1)_{t,H},{\bf
P}_2)_{t,H}$ of twice iterated wrap groups twisted over $ \{
i_0,...,i_{2^r-1} \} $.}

\par {\bf Proof.} Let $M_b$ and $N_b$ be $H^t_p$ manifolds over
${\cal A}_r$ with $1\le r \le 2$, $b=1, 2$ and let $G$ be a group
satisfying Conditions 15$(A1,A2,C1,C2)$ such that $E(N_b,G,\pi ,\Psi
)$ is a principal $G$-bundle. Consider the smashed products
$M_1\otimes _lM_2$, $N_1\otimes _lN_2$ of manifolds and the smashed
product group $G^s$ (see Proposition 19), where $t\ge \max (dim
(M_1), dim (N_1), dim (G)) +1$, where $dim (M_b)$ is a covering
dimension of $M_b$ (see \cite{eng}), $dim (M_1)= dim (M_2)$, $dim
(N_1) = dim (N_2)$. For $At (M_b) = \{ (U_{j,b},\phi _{j,b}): j \} $
an atlas of $M_b$ its connecting mappings $\phi _{j,b}\circ \phi
_{k,b}^{-1}$ are $H^t_p$ functions over ${\cal A}_r$ for
$U_{j,b}\cap U_{k,b}\ne \emptyset $, where $\phi _{j,b}: U_{j,b}\to
{\cal A}_r$ are homeomorphisms of $U_{j,b}$ onto $\phi
_{j,b}(U_{j,b})$. Then $M_1\otimes _lM_2$ consists of all points
$(x,yl)$ with $x\in M_1$ and $y\in M_2$, with the atlas $At
(M_1\otimes _lM_2) = \{ (U_{j,1}\otimes _lU_{q,2}, \phi
_{j,1}\otimes _l\phi _{q,2}): j, q \} $ such that $\phi
_{j,1}\otimes _l \phi _{q,2}: U_{j,1}\otimes _lU_{q,2}\to {\cal
A}_{r+1}^m$, where $m$ is a dimension of $M_1$ over ${\cal A}_r$.
Express for $z=x+yl\in {\cal A}_r$ with $x, y\in {\cal A}_r$ numbers
$x, y$ in the $z$ representation, then denote by $\theta _{j,q}$
mappings corresponding to $\phi _{j,1}\otimes _l\phi _{q,2}$ in the
$z$ representation, hence the transition mappings $\theta
_{j,q}\circ \theta _{k,n}^{-1}$ are $H^t_p$ over ${\cal A}_{r+1}$,
when $(U_{j,1}\otimes _lU_{q,2})\cap (U_{k,1}\otimes _lU_{n,2})\ne
\emptyset $. Therefore, $M_1\otimes _lM_2$ and $N_1\otimes _lN_2$
are $H^t_p$ manifolds over ${\cal A}_{r+1}$.

\par In view of the Sobolev embedding theorem each $H^t$ mapping
on $M_1\otimes _lM_2$ or $N_1\otimes _lN_2$ or $G^s$ is continuous
for $t$ satisfying the inequality \par $t\ge \max (dim (M_1), dim
(N_1), dim (G)) +1$, where $dim (M_1) = dim (M_2)$, $dim (N_1) = dim
(N_2)$.

\par Each locally analytic function $f(x,y) = f_1(x,y) + f_2(x,y)l$ by
$x\in U$ and $y\in V$ can be written as the locally analytic
function by $z=x+yl$ with values in ${\cal A}_{r+1}$, where $U$ and
$V$ are open in ${\cal A}_r^m$, $f_b(x,y)$ is a locally analytic
function with values in ${\cal A}_r^w$, $b=1, 2$, $m, w\in \bf N$.
Indeed, write each variable $x_j$ and $y_j$ through $z_j$ with the
help of generators of ${\cal A}_{r+1}$, where $x_j, y_j\in {\cal
A}_r$, $z_j\in {\cal A}_{r+1}$, $x = (x_1,...,x_m)\in {\cal A}_r^m$,
$z = (z_1,...,z_m)\in {\cal A}_{r+1}^m$ (see Formulas 2.8$(2)$ and
Theorem 2.16 \cite{luoyst2}). If $z\in {\cal A}_{r+1}$, then \par
$(1)$ $z= v_0i_0+...+v_{2^{r+1}-1} i_{2^{r+1}-1}$, where $v_j\in \bf
R$ for each $j=0,...,2^{r+1}-1$, \par $(2)$ $v_0 = (z +
(2^{r+1}-2)^{-1} \{ -z + \sum_{j=1}^{2^{r+1}-1} i_j(zi_j^*) \} )/2$,
\par $(3)$ $v_s = (i_s(2^{r+1}-2)^{-1} \{ - z + \sum_{j=1}^{2^{r+1}-1}
i_j(zi_j^*) \} - z i_j )/2$ for each $s=1,...,2^{r+1}-1$, where $z^*
= {\tilde z}$ denotes the conjugated Cayley-Dickson number $z$. At
the same time we have for $z = x+yl$ with $x, y\in {\cal A}_r$, that
\par $(4)$ $x = v_0i_0+...+v_{2^r-1}i_{2^r-1}$ and \par $(5)$ $y =
(v_{2^r}i_{2^r}+...+v_{2^{r+1}-1}i_{2^{r+1}-1})l^*$, \\ where
$l=i_{2^r}$ denotes the doubling generator.

\par Therefore, $f(x,y)$ becomes
${\cal A}_{r+1}$ holomorphic using the corresponding phrases arising
canonically from expressions of $x_j, y_j$ through $z_j$ by Formulas
$(1-5)$. The set of holomorphic functions is dense in $H^t_p$ in
accordance with the definition of this space, hence using a Cauchy
net we can consider for each $f_1, f_2\in H^t_p$ over ${\cal A}_r$ a
representation of a function $f= f_1+f_2l$ belonging to $H^t_p$ over
${\cal A}_{r+1}$ (see also \cite{luoyst2,lufsqv}).

\par Then $E(N_1\otimes _lN_2,G^s,\pi ^s,\Psi ^s)$ is naturally isomorphic
with $E(N_1,G,\pi _1,\Psi _1)\otimes _l E(N_2,G,\pi _2,\Psi _2)$,
where $\pi ^s =\pi _1\otimes \pi _2l: E(N_1\otimes _lN_2,G^s,\pi
^s,\Psi ^s)\to N_1\otimes _l N_2$ is the natural projection.

\par If $\gamma : M_1\otimes _lM_2\to N_1\otimes _lN_2$ is an
$H^t_p$ mapping, then $\gamma (z) = \gamma _1(x,y) \times \gamma
_2(x,y)l$, where $x\in M_1$ and $y\in M_2$, $z=(x,yl)\in M_1\otimes
_lM_2$, $\gamma _b: M_1\otimes _lM_2\to N_b$. We can write $\gamma
_b(x,y)$ as $(\gamma _{b,1}(x))(y)$ a family of functions by $x$ and
a parameter $y$ or as $(\gamma _{b,2}(y))(x)$ a family of functions
by $y$ with a parameter $x$. If $\eta _{b,a}: M_a\to N_b$, then
${\bf P}_{{\hat {\eta }}_{b,a},u_b,a}$ denotes the parallel
transport structure on $M_a$ over $E(N_b,G,\pi _b,\Psi _b)$.

\par Then ${\bf P}^s_{{\hat {\gamma }},u}
(z) = [{\bf P}_{{\hat {\gamma }}_{1,1},u_1;1}(x)] [{\bf P}_{{\hat
{\gamma }}_{1,2},u_1;2}(y)] \otimes _l [{\bf P}_{{\hat {\gamma
}}_{2,1},u_2;2}(x)] [{\bf P}_{{\hat {\gamma }}_{2,2},u_2;2}(y)] \in
E_{y_0}(N_1\otimes _lN_2,G^s,\pi ^s, \Psi ^s)$ is the parallel
transport structure in $M_1\otimes _lM_2$ induced by that of in
$M_1$ and $M_2$, where $u\in E_{y_0} (N_1\otimes _lN_2, G^s, \pi ^s,
\Psi ^s)$, $u= u_1\otimes _l u_2$, $u_b\in E_{y_{0,b}}(N_b,G,\pi
_b,\Psi _b)$, $y_{0,b}\in N_b$ is a marked point, $b=1, 2$, $y_0 =
y_{0,1}\otimes _l y_{0,2}$. Then ${\bf P}^s$ is $G^s$ equivariant.
Therefore, $<{\bf P}^s_{{\hat {\gamma }},u}>_{t,H} = < {\bf
P}_{{\hat {\gamma }}_1,u_1}>_{t,H} \otimes _l <{\bf P}_{{\hat
{\gamma }}_2,u_2}>_{t,H} = <[{\bf P}_{{\hat {\gamma
}}_{1,1},u_1;1}(x)] [{\bf P}_{{\hat {\gamma
}}_{1,2},u_1;2}(y)]>_{t,H} \otimes _l <[{\bf P}_{{\hat {\gamma
}}_{2,1},u_2;2}(x)] [{\bf P}_{{\hat {\gamma
}}_{2,2},u_2;2}(y)]>_{t,H}$, where ${\bf P}_{{\hat {\gamma
}}_b,u_b}$ is the parallel transport structure on $M_1\otimes _lM_2$
over $E(N_b,G,\pi _b,\Psi _b)$, $b=1, 2$.

\par Hence $(W^{M_1\otimes _lM_2; \{
s_{0,j,1}\otimes _l s_{0,v,2}: j=1,...,k_1; v=1,...,k_2 \} }
E;N_1\otimes _lN_2,G^s,{\bf P}^s)_{t,H}$ is isomorphic with the
smashed product \par $W^{M_2; \{ s_{0,v,2}: v=1,...,k_2 \}
}E;N_1,(W^{M_1; \{ s_{0,j,1}: j=1,...,k_1 \} }E;N_1,G,{\bf
P}_1)_{t,H},{\bf P}_2)_{t,H}\otimes _l$ \par $ W^{M_2; \{ s_{0,v,2}:
v=1,...,k_2 \} }E;N_2,(W^{M_1; \{ s_{0,j,1}: j=1,...,k_1 \}
}E;N_2,G,{\bf P}_1)_{t,H},{\bf P}_2)_{t,H}$ of iterated wrap groups.

\par {\bf 21. Theorem.} {\it There exists a homomorphism of iterated
wrap groups $\theta : (W^ME)_{a;\infty ,H}\otimes (W^ME)_{b;\infty
,H}\to (W^ME)_{a+b;\infty ,H}$ for each $a, b\in \bf N$, where $G$
is an $H^{\infty }_p$ group, $E(N,G,\pi ,\Psi )$ is the principal
$H^{\infty }_p$ bundle with the structure group $G$. Moreover, if
$G$ is either associative or alternative, then $\theta $ is either
associative or alternative.}

\par {\bf Proof.} Consider iterated wrap groups $(W^ME)_{a;\infty ,H}$
as in \S 4, $a\in \bf N$. If $\gamma _a: M^a\to N$, $\gamma _b:
M^b\to N$ are $H^{\infty }_p$ mappings such that $\gamma
_b(s_{0,j_1}\times ...\times s_{0,j_b})=y_0$ for each $j_l=1,...,k$
and $l=1,...,b$, then $\gamma := \gamma _a\times \gamma _b:
M^a\times M^b\to N\times N=N^2$, where $M^a\times M^b = M^{a+b}$,
$s_{0,j}$ are marked points in $M$ with $j=1,...,k$ and $y_0$ is a
marked point in $N$, $H^{\infty }_p = \bigcap_{t\in \bf N}H^t_p$.
This gives the iterated parallel transport structure ${\bf P}_{{\hat
{\gamma }},u;a+b}(x) := {\bf P}_{{\hat {\gamma
}}_a,u_a;a}(x_a)\otimes {\bf P}_{{\hat {\gamma }}_b,u;b}(x_b)$ on
$M^{a+b}$ over $E(N^2,G^2,\pi ,\Psi )$, where $u_b\in
E_{y_0}(N,G,\pi ,\Psi )$, $u = u_a\times u_b \in E_{y_0\times y_0}
(N^2,G^2,\pi ,\Psi )$.
\par The bunch $M^b\vee M^b$ is taken by points $s_{j_1,...,j_b}$
in $M^b$, where $s_{j_1,...,j_b} := s_{0,j_1}\times ... \times
s_{0,j_b}$ with $j_1,...,j_b\in \{ 1,...,k \} $; $s_{0,j}$ are
marked points in $M$ with $j=1,...,k$. Then $(M^a\vee M^a)\times
(M^b\vee M^b)\setminus \{ s_{j_1,...,j_{a+b}}: j_l=1,...,k;
l=1,...,a+b \} $ is $H^t_p$ homeomorphic with $M^{a+b}\vee
M^{a+b}\setminus \{ s_{j_1,...,j_{a+b}}: j_l=1,...,k; l=1,...,a+b \}
$, since $s_{j_1,...,j_a}\times s_{j_{a+1},...,j_{a+b}} =
s_{j_1,...,j_{a+b}}$ for each $j_1,...,j_{a+b}$. There is the
embedding $DifH^{\infty }_p(M^a)\times DifH^{\infty
}_p(M^b)\hookrightarrow DifH^{\infty }_p(M^{a+b})$ for each $a, b\in
\bf N$. If $f_a \in DifH^{\infty }_p(M^a)$ having a restriction
$f_a|_{K_a}=id$, then $f_a\times f_b\in DifH^{\infty }_p(M^{a+b})$
and $f_a\times f_b|_{K_a\times K_b}=id$ for $K_a\subset M^a$. Put
$\theta ( <{\bf P}_{{\hat {\gamma }}_a,u_a;a}>_{\infty ,H;a}, <{\bf
P}_{{\hat {\gamma }}_b,u_b;b}>_{\infty ,H;b}) = < <{\bf P}_{{\hat
{\gamma }}_a,u_a;a}>_{\infty ,H;a}\otimes <{\bf P}_{{\hat {\gamma
}}_b,u_b;b}>_{\infty ,H;b}>_{\infty ,H;a+b}$ is the group
homomorphism, where the detailed notation $<*>_{t,H;a}$ denotes the
equivalence class over the manifold $M^a$ instead of $M$, $a\in \bf
N$.

\par Therefore,
$<{\bf P}_{{\hat {\gamma }}\vee {\hat {\eta }},u;a+b}>_{\infty
,H;a+b} := < <{\bf P}_{{\hat {\gamma }}_a\vee {\hat {\eta
}}_a,u_a;a}>_{\infty ,H;a}\otimes <{\bf P}_{{\hat {\gamma }}_b\vee
{\hat {\eta }}_b,u_b;b}>_{\infty ,H;b}>_{\infty ,H;a+b} $ \\
$=< (<{\bf P}_{{\hat {\gamma }}_a,u_a;a}>_{\infty ,H;a} <{\bf
P}_{{\hat {\eta }}_a,u_a;a}>_{\infty ,H;a})\otimes  (<{\bf P}_{{\hat
{\gamma }}_b,u_b;b}>_{\infty ,H;b} <{\bf P}_{{\hat {\eta
}}_b,u_b;b}>_{\infty ,H;b})>_{\infty ,H;a+b}$ \\
$=< (<{\bf P}_{{\hat {\gamma }}_a,u_a;a}>_{\infty ,H;a}\otimes <{\bf
P}_{{\hat {\gamma }}_b,u_b;b}>_{\infty ,H;b})  (<{\bf P}_{{\hat
{\eta }}_a,u_a;a}>_{\infty ,H;a}\otimes <{\bf P}_{{\hat {\eta
}}_b,u_b;b}>_{\infty ,H;b})>_{\infty ,H;a+b}$ \\
$=< <{\bf P}_{{\hat {\gamma }}_a,u_a;a}>_{\infty ,H;a}\otimes <{\bf
P}_{{\hat {\gamma }}_b,u_b;b}>_{\infty ,H;b}>_{\infty ,H;a+b} <
<{\bf P}_{{\hat {\eta }}_a,u_a;a}>_{\infty ,H;a}\otimes <{\bf
P}_{{\hat {\eta
}}_b,u_b;b}>_{\infty ,H;b} >_{\infty ,H;a+b}$ \\
$=\theta ( <{\bf P}_{{\hat {\gamma }}_a,u_a;a}>_{\infty ,H;a}, <{\bf
P}_{{\hat {\gamma }}_b,u_b;b}>_{\infty ,H;b}) \theta ( <{\bf
P}_{{\hat {\eta }}_a,u_a;a}>_{\infty ,H;a}, <{\bf P}_{{\hat {\eta
}}_b,u_b;b}>_{\infty ,H;b})$. \\
Thus $\theta $ is the group homomorphism.
\par The mapping $H^{\infty }_p(M^a,N)\times H^{\infty }_p(M^b,N)\ni
(\gamma _a\times \gamma _b)\mapsto (\gamma _a,\gamma _b)\in
H^{\infty }_p(M^{a+b},N^2)$ is of $H^{\infty }_p$ class. The
multiplication in $G^v$ is $H^{\infty }_p$ for each $v\in \bf N$,
since it is such in $G$, since the multiplication in $G^v$ is
$(a_1,...,a_v)\times (b_1,...,b_v)=(a_1b_1,...,a_vb_v)$, where $G^v$
is the $v$ times direct product of $G$, $a_1,...,a_v, b_1,...,b_v\in
G$.
\par The iterated wrap group $(W^ME)_{l;t,H}$ for the bundle $E$ is
the principal $G^{kl}$ bundle over the iterated commutative wrap
group $(W^MN)_{l;t,H}$ for the manifold $N$, since the number of
marked points in $M^l$ is $kl$, where $E$ is the principal $G$
bundle on the manifold $N$, $l\in \bf N$. Thus the iterated wrap
group is associative or alternative if such is $G$. In view of
Proposition 7 and Remark 4 the homomorphism $\theta $ is of
$H^{\infty }_p$ class. From the wrap monoids it has the natural
$H^{\infty }_p$ extension on wrap groups.
\par If $G$ is associative, then \par $<{\bf P}_{{\hat {\gamma
}},u;a+b+v}>_{\infty ,H;a+b+v} = < < ( <{\bf P}_{{\hat {\gamma
}}_a,u_a; a}>_{\infty ,H;a} \otimes <{\bf P}_{{\hat {\gamma
}}_b,u_b;b}>_{\infty ,H;b})
>_{\infty ,H;a+b} \otimes <{\bf P}_{{\hat {\gamma }}_v,u_v; v}
>_{\infty ,H;v}>_{\infty ,H;a+b+v}$ \\ $ = <
<{\bf P}_{{\hat {\gamma }}_a,u_a;a}>_{\infty ,H;a}\otimes ( <{\bf
P}_{{\hat {\gamma }}_b,u_b;b}>_{\infty ,H;b}\otimes <{\bf P}_{{\hat
{\gamma }}_v,u_v; v}>_{\infty ,H;v}) >_{\infty ,H;a+b+v} = \theta (
\theta ( <{\bf P}_{{\hat {\gamma }}_a,u_a; a}>_{t,H;a}, <{\bf
P}_{{\hat {\gamma }}_b,u_b;b}>_{t,H;b}), <{\bf P}_{{\hat {\gamma
}}_v,u_v; v}>_{t,H;v}) $  \\
$\theta ( <{\bf P}_{{\hat {\gamma }}_a,u_a; a}>_{t,H;a}, \theta (
<{\bf P}_{{\hat {\gamma }}_b,u_b;b}>_{t,H;b}), <{\bf P}_{{\hat
{\gamma }}_v,u_v; v}>_{t,H;v})) $, \\
consequently, $\theta $ is the associative homomorphism.
\par If $G$ is alternative, then
\par $<{\bf P}_{{\hat {\gamma }},u;a+a+b}>_{\infty ,H;a+a+b} = < < (
<{\bf P}_{{\hat {\gamma }}_a,u_a; a}>_{\infty ,H;a} \otimes <{\bf
P}_{{\hat {\gamma }}_a,u_a;a}>_{\infty ,H;a})
>_{\infty ,H;a+a} \otimes <{\bf P}_{{\hat {\gamma }}_b,u_b; b}
>_{\infty ,H;v}>_{\infty ,H;a+a+b}$ \\ $ = <
<{\bf P}_{{\hat {\gamma }}_a,u_a;a}>_{\infty ,H;a}\otimes ( <{\bf
P}_{{\hat {\gamma }}_a,u_a;a}>_{\infty ,H;a}\otimes <{\bf P}_{{\hat
{\gamma }}_b,u_b; b}>_{\infty ,H;b}) >_{\infty ,H;a+a+b} = \theta (
\theta ( <{\bf P}_{{\hat {\gamma }}_a,u_a; a}>_{t,H;a}, <{\bf
P}_{{\hat {\gamma }}_a,u_a;a}>_{t,H;a}), <{\bf P}_{{\hat {\gamma
}}_b,u_b; b}>_{t,H;b}) $  \\
$\theta ( <{\bf P}_{{\hat {\gamma }}_a,u_a; a}>_{t,H;a}, \theta (
<{\bf P}_{{\hat {\gamma }}_a,u_a;a}>_{t,H;a}), <{\bf P}_{{\hat
{\gamma }}_b,u_b; b}>_{t,H;b})) $, \\
consequently, the homomorphism $\theta $ is alternative from the
left, analogously it is alternative from the right.


\begin{thebibliography}{99}

\bibitem{bredon} G.E. Bredon. "Sheaf theory" (New York: McGraw-Hill,
1967).

\bibitem{dubnovfom} B.A. Dubrovin, S.P. Novikov, A.T. Fomenko.
"Modern geometry" (Moscow: Nauka, 1979).

\bibitem{emch} G. Emch. Helv. Phys. Acta.
"M$\grave e$chanique quantique quaternionienne et
Relativit$\grave e$ restreinte", {\bf 36} (1963), 739-788.

\bibitem{eng} R. Engelking. "General topology" (Moscow: Mir, 1986).

\bibitem{gaj} P. Gajer. "Higher Holonomies, Geometric Loop Groups and Smooth
Deligne Cohomology". in: "Advances in Geometry". J.-L. Brylinski ed.
Progr. Math. V. {\bf 172}, P. 195-235 (Boston: Birkha\"user, 1999).

\bibitem{guetze} F. G\"ursey, C.-H. Tze. "On the role of
division, Jordan and related algebras in particle physics"
(Singapore: World Scientific Publ. Co., 1996).

\bibitem{hamilt} W.R. Hamilton. "Selected papers. Optics. Dynamics.
Quaternions" (Moscow: Nauka, 1994).

\bibitem{harvey} F.R. Harvey. "Spinors and calibrations".
Perspectives in Mathem. {\bf 9} (Boston: Academic Press, 1990).

\bibitem{ish} C.J. Isham. "Topological and global aspects of
quantum theory". In: "Relativity, groups and topology.II" 1059-1290,
(Les Hauches, 1983). Editors: R. Stora, B.S. De Witt (Amsterdam:
Elsevier Sci. Publ., 1984).

\bibitem{kansol} I.L. Kantor, A.S. Solodovnikov.
"Hypercomplex numbers" (Berlin: Springer-Verlag, 1989).

\bibitem{lawmich} H.B. Lawson, M.-L. Michelson. "Spin geometry"
(Princeton: Princ. Univ. Press, 1989).

\bibitem{luwrgfbqo} S.V. Ludkovsky.
"Wrap groups of fiber bundles over quaternions and octonions". Los
Alam. Nat. Lab. {\bf math.FA 0802.0661}, 27 pages.

\bibitem{ludan} S.V. Ludkovsky. Dokl. Akad. Nauk.
"Quasi-invariant measures on loop groups of Riemann manifolds",
{\bf 370: 3} (2000), 306-308.

\bibitem{lupom} S.V. Ludkovsky. Southeast Asian Bulletin of
Mathematics. "Poisson measures for topological groups and their
representations", {\bf 25} (2002), 653-680. (shortly in Russ. Math.
Surv. {\bf 56: 1} (2001), 169-170; previous versions: {\bf
IHES/M/98/88}, 38 pages, also Los Alamos Nat. Lab. {\bf
math.RT/9910110}).

\bibitem{lufsqv} S.V. Ludkovsky. J. Mathem. Sci.
"Functions of several Cayley-Dickson variables and manifolds over
them", {\bf 141: 3} (2007), 1299-1330 (previous variant: Los Alamos
Nat. Lab. {\bf math.CV/0302011}).

\bibitem{norfamlud} S.V. Ludkovsky.  Sovrem.
Mathem. Fundam. Napravl. "Normal families of functions and groups of
pseudoconformal diffeomorphisms of quaternion and octonion
variables", {\bf 18} (2006), 101-164 (previous variant: Los Alam.
Nat. Lab. math.DG/0603006).

\bibitem{luoyst} S.V. Ludkovsky, F. van Oystaeyen. Bull. Sci. Math.
(Paris). Ser. 2. "Differentiable functions of quaternion variables",
{\bf 127} (2003), 755-796.

\bibitem{luoyst2} S.V. Ludkovsky.  J. Mathem. Sci.
"Differentiable functions of Cayley-Dickson numbers and line
integration", {\bf 141: 3} (2007), 1231-1298 (previous version: Los
Alam. Nat. Lab. math.NT/0406048; math.CV/0406306; math.CV/0405471).

\bibitem{lujmslg} S.V. Ludkovsky. J. Mathem. Sci. "Stochastic
processes on geometric loop groups, diffeomorphism groups of
connected manifolds, associated unitary representations", {\bf 141:
3} (2007), 1331-1384 (previous version: Los Alam. Nat. Lab.
math.AG/0407439,  July 2004).

\bibitem{lufoclg} S.V. Ludkovsky.
"Geometric loop groups and diffeomorphism groups of manifolds,
stochastic processes on them, associated unitary representations".
In the book: "Focus on Groups Theory Research" (Nova Science
Publishers, Inc.: New York) 2006, pages 59-136.

\bibitem{lugmlg} S.V. Ludkovsky. J. Mathem. Sci.
"Generalized geometric loop groups of complex manifolds, Gaussian
quasi-invariant measures on them and their representations", {\bf
122: 1} (2004), 2984-3011 (earlier version: Los Alam. Nat. Lab. {\bf
math.RT/9910086}, October 1999).

\bibitem{lufejms} S.V. Ludkovsky. Far East J. of Math.
Sci. (FJMS). "Quasi-conformal functions of quaternion and octonion
variables, their integral transformations", {\bf 28: 1} (2008),
37-88.

\bibitem{mensk} M.B. Mensky. "The paths group. Measurement. Fields.
Particles" (Moscow: Nauka, 1983).

\bibitem{michor} P.W.  Michor.  "Manifolds of Differentiable
Mappings" (Boston: Shiva, 1980).

\bibitem{miha} V.P. Mihailov. "Differential equations in
partial derivatives" (Moscow: Nauka, 1976).

\bibitem{milmorse} J. Milnor. "Morse theory" (Princeton, New Jersey:
Princeton Univ. Press, 1963).

\bibitem{omo} H. Omori. Trans. Amer. Math. Soc.
"Groups of diffeomorphisms and their subgroups", {\bf 179} (1973),
85-122.

\bibitem{omori2} H. Omori. J. Math. Soc. Japan.
"Local structures of groups of diffeomorphisms",  {\bf 24: 1}
(1972), 60-88.

\bibitem{pont} L.S. Pontrjagin. "Continuous groups"
(Moscow: Nauka, 1984).

\bibitem{seel} R.T. Seeley. Proceed.  Amer.  Math.  Soc.
"Extensions of $C^{\infty }$ Functions Defined in a Half Space",
{\bf 15} (1964), 625-626.

\bibitem{souriau} J.M. Souriau. "Groupes differentiels"
(Berlin: Springer Verlag, 1981).

\bibitem{steen} N. Steenrod. "The topology of fibre budles"
(Princeton, New Jersey: Princeton Univ. Press, 1951).

\bibitem{sulwint} R. Sulanke, P. Wintgen. "Differentialgeometrie
und Faserb\"undel" (Berlin: Veb deutscher Verlag der Wissenschaften,
1972).

\bibitem{swit} R.M.  Switzer.  "Algebraic Topology - Homotopy and Homology"
(Berlin: Springer-Verlag, 1975).

\bibitem{touger} J.C.  Tougeron.  "Ideaux de
Fonctions Differentiables" (Berlin: Springer-Verlag, 1972).

\end{thebibliography}
\end{document}